\documentclass[10pt]{article}
\usepackage{amsmath,stackengine, leftidx}
\usepackage{amssymb}
\usepackage{titlesec}
\usepackage[T1]{fontenc}
\usepackage{empheq}
\usepackage{graphicx}
\usepackage{titling}
\usepackage[english]{babel}
\usepackage[margin=2.5cm]{geometry}
\usepackage{csquotes}
\usepackage{lipsum}
\usepackage[english]{babel}
\usepackage{amsthm}
\usepackage{amssymb}
\usepackage{bm} 
\usepackage{hyperref}
\usepackage{biblatex} 
\usepackage{romannum}
\usepackage{mathrsfs}
\usepackage{color}
\usepackage{mathtools,amsfonts,amssymb,setspace}
\usepackage{enumitem}
\usepackage{subfigure}
\usepackage[font=small,labelfont=bf]{caption}
\newtheorem{theorem}{Theorem}[section]
\newtheorem{corollary}{Corollary}[theorem]
\newtheorem{lemma}[theorem]{Lemma}
\theoremstyle{definition}
\newtheorem{definition}{Definition}[section]
\newtheorem{prop}[theorem]{Proposition}
\theoremstyle{definition}
\newtheorem{coro}[corollary]{Corollary}
\theoremstyle{definition}
\newtheorem{remark}[theorem]{Remark}

\newcommand{\subjclass}[1]{\small	\par\noindent\textbf{Mathematics Subject Classification (2020): }#1\par}
\newcommand{\keywords}[1]{\small	\par\noindent\textbf{Keywords: }#1\par}
\definecolor{green}{rgb}{0.1,0.62,0.0}
\definecolor{blue}{rgb}{0.0, 0.0, 1.0}
\definecolor{deepblue}{rgb}{0.0,0.0,0.7}
\definecolor{red}{rgb}{1.0, 0.0, 0.0}

\font\myfont=cmbx12 at 11pt
\font\myfond=cmr12 at 10pt
\title{\myfont CONTROLLABILITY OF THE CAHN-HILLIARD SYSTEM AT FINITELY MANY OBSERVATION POINTS  }
\author{\myfond SHEETAL DHARMATTI \thanks{sheetal@iisertvm.ac.in \; School of Mathematics, Indian Institute of Science Education and Research, Thiruvananthapuram, Maruthamala PO, Vithura, Thiruvananthapuram, Kerala, 695551, INDIA}  \;  and  \;  GREESHMA K \thanks{ greeshmak21@iisertvm.ac.in \; School of Mathematics, Indian Institute of Science Education and Research, Thiruvananthapuram, Maruthamala PO, Vithura, Thiruvananthapuram, Kerala, 695551, INDIA,\\ \textbf{Acknowledgement}: Greeshma K would like to thank the Department of Science and Technology (DST), India, for the Innovation in Science Pursuit for Inspired Research (INSPIRE) Fellowship (IF210199). \\
\textbf{Declaration}: Authors declare that there is no conflict of interest. Both authors have contributed equally in designing and resolving of the problem and preparation of the manuscript.} }

\date{}

\numberwithin{equation}{section}
\thanksmarkseries{arabic}

\begin{document}

\maketitle
\begin{abstract}
We study approximate controllability of the  Cahn–Hilliard (CH) system with prescribed targets at finitely many points in the spatial domain. Our aim is to design a $L^2$ control such that given a finite set of observation points and target values, the solution reaches arbitrarily close to the prescribed targets at these points. First we  study the associated linearised system using variational methods and then extend the results to the semi-linear CH system using fixed-point arguments. Further, we investigate approximate controllability in more general setting where we allow the control to be measure valued (in specific, Dirac delta type controls). This formulation is more challenging  to obtain requisite regularity due to the simultaneous presence of pointwise observations and pointwise controls. Hence we restrict ourselves to the study of linearised CH system and devise a control which exhibits approximate controllability. We achieve this by approximating measure valued control with $L^2$ controls and showing appropriate convergence.

\end{abstract}

\pagenumbering{arabic}
\setcounter{section}{1}
\section*{Introduction}

Let $\Omega$ be an open, bounded subset of $\mathbb{R}^2$ with a smooth boundary $\partial\Omega$. Let $O$ be an open subset of $\Omega$ and denote $O \times [0,T]$ by $q$. Consider the Cahn-Hilliard system given by, 
\begin{align}
    \varphi' + (\Delta)^2\varphi - \Delta F'(\varphi) &= U\chi_{q}, \text{ in } \Omega \times(0,T),\label{eq91a}\\
    \frac{\partial\varphi}{\partial\textbf{n}}=  \frac{\partial\Delta\varphi}{\partial\textbf{n}} &= 0, \text{ on } \partial\Omega\times(0,T),\label{eq91b}\\ 
    \varphi(0) &= \varphi_0, \text{ in } \Omega.\label{eq91c}
\end{align}
Here, $\varphi$ denotes the order parameter and $F$ is a regular potential of polynomially controlled growth. The double well potential given by, $F(s) = (1-s^2)^2/2$ is an example for a regular potential. In the present work, we restrict ourselves to a class of potentials whose derivative is sublinear at infinity. The double-well potential does not satisfy this assumption since $F'$ has cubic growth. On the other hand, the quadratic potential $F(s)=\frac{1}{2}(1-s^2)$ satisfies the required linear growth condition.

We study a pointwise approximate controllability problem for the semilinear Cahn-Hilliard equation. More precisely, we seek controls that steer the solution arbitrarily close to prescribed target values at finitely many points in the spatial domain. Let $x_1,x_2,\ldots,x_k\in\Omega$ be a finite collection of observation points, and let
$\omega_1, \omega_2,..., \omega_k \in \mathbb{R}$ denotes respective desired state at time $T$. Given $\epsilon >0$, our aim is to find a control $U\in L^2(q)$ such that the solution of the Cahn-Hilliard system \eqref{eq91a}- \eqref{eq91c} satisfies,
\begin{equation*}
    \sum_{i=1}^{k}(\varphi(x_i,T)-\omega_i)^2  \leq \epsilon. 
\end{equation*}

Pointwise observations and controls are motivated by a variety of practical applications in which measurements or control inputs are available only at a finite number of spatial locations. In such settings, pointwise strategies may be more practical and cost-effective than distributed observations or controls. Moreover, by appropriately selecting the observation points, one may obtain accurate tracking at the locations of interest without requiring the solution to remain close to the desired state throughout the entire domain.

Our work is motivated by the approximate controllability problem for semilinear heat equations studied in \cite{Fabre95}. In an approximate controllability problem, the objective is to find a control that steers the solution arbitrarily close to a prescribed target state. In \cite{Fabre95}, the authors considered two control problems: one in which both the control and the target state belong to $L^2(\Omega)$, and another in which the control belongs to $L^2(\Omega)$, while the observations are made in $C(\Bar{\Omega})$. Their results were established under the condition that the nonlinearity is globally Lipschitz. Their analysis is based on a variational approach, in which the approximate controllability problem is reformulated as the minimisation of a suitable functional involving the adjoint variable.

Similar ideas were employed in \cite{DRA} to study approximate controllability for higher order Cahn-Hilliard type equations. The unique continuation property of the adjoint system obtained from \cite{saut87} is crucial in the analysis. We adopt a similar variational framework in the present work. Since the observations are prescribed only at finitely many points of the domain, the adjoint system involves source terms of lower regularity, which introduces additional analytical challenges. In contrast to \cite{DRA}, where the desired state belongs to $L^2(\Omega)$, our desired state is specified through its values at finitely many points. This pointwise formulation constitutes the novelty of the present work. A related approximate controllability problem for semilinear heat equations was considered in \cite{FCZ}, where the growth assumptions on the nonlinearity were slightly weakened. As in \cite{Fabre95}, the analysis is based on reformulating the controllability problem as a minimisation problem in the dual space.

The present work is also motivated by the sidewise tracking controllability problem for the heat equation studied in \cite{EZP}. In \cite{EZP}, authors have considered approximate tracking controllability with boundary controls and observations restricted to a portion of the boundary. Thus, both the control and observation sets are of measure zero. By combining variational techniques with unique continuation arguments, they established approximate tracking results for the corresponding problem. Although the setting in our problem is different, the approach is inspired by the ideas in \cite{EZP}. We also mention the recent work \cite{ABM}, where a related tracking controllability problem was studied for semilinear parabolic equations. In the one-dimensional setting, the desired trajectory is prescribed at finitely many spatial points as a function of time. The authors studied both constant target states and time-dependent target trajectories (moving target). 

Regarding the Cahn-Hilliard equation, the exact, null, and approximate controllability results are available under various restrictive assumptions on the potential. In \cite{GAO,KAS,KAK}, the authors derived Carleman estimates for fourth-order semilinear parabolic equations and subsequently obtained controllability results for the corresponding systems. In the present work, instead of employing a Carleman estimate, we make use of the unique continuation result established in \cite{saut87}.

We consider two control problems. In the first part of the paper, we study approximate controllability with pointwise observations and distributed controls $U \in L^2(q)$. We establish the approximate controllability of the corresponding linearised Cahn-Hilliard system by variational methods where the idea is originally from \cite{ LIO, LIO1990}. Further we extend the result to the semilinear equation by applying Kakutani's fixed point theorem.

In the second part, we consider controls acting at finitely many spatial points and investigate the corresponding approximate controllability problem with pointwise observations for the linearised system. The latter formulation leads to a system with low regular source terms, and obtaining higher regularity of the solution is even more difficult. The main challenges of our work are the following.

\begin{itemize}
    \item A higher regularity of the adjoint system with a low regular terminal data is necessary to make use of the unique continuation property of the adjoint system from \cite{saut87}.
    \item A higher regularity of the linearised variable associated to a measure source data is essential to define the controllability problem. 
    \item Obtaining a pointwise control which solves the approximate controllability problem  as a limit of $L^2(q)$ controls which solves a specific $L^2$ approximate controllability problem is less trivial. 
\end{itemize}

The paper is organised as follows. In the following section, we introduce the notations and functional framework to approximate controllability of the CHB system  with pointwise observation at terminal time $T$. The later part of this section is devoted to the study of approximate controllability problem with distributed $L^2$ controls. For, we first analyse the linearised system and subsequently establish the corresponding result for the semilinear equation using fixed point argument. In Section 3, we investigate CHB system with pointwise observations at terminal time  and  pointwise control acting at finitely many spatial points. However, we have been successful to obtain the results for the linearised CHB system due to lower regular control.
\\\\\textbf{Preliminaries:} Let $\Omega$ be an open, bounded, connected subset of $\mathbb{R}^2$ with a smooth boundary. We denote the $L^2(\Omega)$ space by $H$, and its norm and inner product by $\|.\|$, $(., .)$ respectively. Further, $V$ represents the $H^1(\Omega)$ space: $\|.\|_{V}$ and $(.,.)$ denote the corresponding norm and inner product in the space. We define,
\begin{align*}
     W_N^{2,p}(\Omega) = \{ f\in W^{2,p}(\Omega): \frac{\partial f}{\partial \textbf{n}} =  \frac{\partial \Delta f}{\partial \textbf{n}} = 0 \}. 
\end{align*}
and $W_{N}^{-2,p'}(\Omega)$ denotes the corresponding dual space. Let $X, Y, Z$ be Banach spaces such that $X, Y \subseteq Z$. We define a Banach space $X+Y$ by,
\begin{equation*}
    X+ Y = \{ x+y : x\in X \text{ and } y\in Y \},
\end{equation*}
with norm, 
\begin{align*}
    \|z\|_{X+Y} := \underset{\substack{z= x+y, \\ x\in X, y\in Y}}{\inf} \{ \|x\|_{X} +\|y\|_{Y} \}.
\end{align*}
\textit{Interpolation space: } Let $0< \theta < 1$, $1\leq p < \infty$, we define the interpolation space $(X, Y)_{\theta, p}$, the collection of elements $z\in X+Y$ such that the following norm is finite. 
\begin{align*}
    \|z\|_{(X,Y)_{\theta, p}} := \Big(\int_{0}^{\infty} K(t,z)^p t^{-(1+\theta p)} dt \Big)^{1/p}, 
\end{align*}
where,
\begin{align*}
    K(t,z) := \underset{\substack{z= x+y, \\ x\in X, y\in Y}}{\inf} \|x\|_X + t\|y\|_Y.
\end{align*}
The above definition is based on Petree's K-method.
\\\textit{Besov space: } Let $0< s< \infty$, $1\leq p < \infty$, $1\leq q\leq \infty$ and $m$ be the smallest integer greater than $s$. Then the space $B^s_{p,q}$ is an intermediate space between $L^p(\Omega)$ and $W^{m,p}(\Omega)$ defined by,
\begin{equation*}
    B^{s}_{p,q} := (L^p(\Omega), W^{m,p}(\Omega))_{s/m, q}.
\end{equation*}

We denote the set of $C^\infty$ functions in $\Omega$ with compact support by $\mathcal{D}(\Omega)$. The set $C^m_b(\Omega)$ consists of functions $f\in C^m(\Omega)$ such that $D^\alpha(f)$ is bounded for $0\leq |\alpha| \leq m$. In particular, $C_b(\Omega)$ denotes the collection of bounded continuous functions on $\Omega$. Further, the space $\mathcal{M}(\Omega) = C(\Bar{\Omega})'$ can be identified with the Banach space of real and Borel measurable function with the total variation norm.

\section{Approximate controllability problem with a distributed \texorpdfstring{$L^2$}{e} control}
\textbf{Problem statement}: Consider the local Cahn-Hilliard system given in \eqref{eq91a}-\eqref{eq91c}. Let $ \mathcal{D} = \{ x_1, x_2,..., x_k\}$ be a finite subset of $\Omega$ and $q= O \times [0,T]$, where $O$ is an open subset of $\Omega$. Given $\epsilon >0$ and $\omega = (\omega_1, ..., \omega_k) \in \mathbb{R}^k$, $k\in \mathbb{N}$, our aim is to find a control $U \in L^2(q)$ such that the solution of the system \eqref{eq91a}- \eqref{eq91c} satisfies,
\begin{equation}\label{eq092}
    \sum_{i=1}^{k}(\varphi(x_i,T)-\omega_i)^2 \leq \epsilon. 
\end{equation}

Before addressing the approximate controllability, we discuss the wellposedness of the system \eqref{eq91a}- \eqref{eq91c}. In literature, the well-posedness of local Cahn-Hilliard system is studied under various assumptions on the potential. Let $\varphi_0\in H^2(\Omega)$ and $U\in L^2(q)$, then there exists a solution of the system \eqref{eq91a}-\eqref{eq91c} with the following regularity (See \cite[Chapter 3]{AML}).
$$\varphi \in L^2(0,T;H^4(\Omega)) \cap L^\infty(0,T;H^2(\Omega)) \cap H^1(0,T;H). $$
Moreover, by interpolation embedding we have $\varphi \in L^2(0,T;H^4(\Omega)) \cap H^1(0,T; H)) \hookrightarrow C([0,T]; H^2(\Omega)),$ further by Sobolev embedding, we have $H^2(\Omega) \hookrightarrow C(\Bar{\Omega})$. Therefore, the pointwise evaluation of $\varphi$ in \eqref{eq092} is well defined. Here we additionally assume that the potential $F \in C^1(\mathbb{R})$ is sublinear at infinity. Precisely,
\begin{enumerate}[font={\bfseries},label={A\arabic*.}]
    \item[{[A1]}] The potential $F \in C^1(\mathbb{R})$ and $\exists C_1, C_2 >0$ such that 
    \begin{equation*}
        F'(s) \leq C_1|s| + C_2, \,\, \forall s\in \mathbb{R}.
    \end{equation*}
    Furthermore, the potential $F$ is twice differentiable at $0$ and $F'(0) = 0$.
\end{enumerate}

\subsection{The linear problem}

Our idea is to linearise the system, study the associated control problem, and then extend the results to the nonlinear case by applying fixed point arguments. With this motive, we introduce a function $g$ as follows. 
\begin{align*}
    g(s) = \begin{cases} 
             \frac{F'(s)}{s}, & \text{ if }s\neq 0,\\
             F''(0), & \text{ if }s=0. 
             \end{cases} 
\end{align*}
Observe that, by assumption $\textbf{[A1]}$, $\|g(s)\|_{L^\infty} \leq C$ for every $s\in \mathbb{R}$; Moreover, $g$ is continuous. We can write,
\begin{equation*}
    F'(s) =  g(s)s, \,\,\forall s\in \mathbb{R}.
\end{equation*}
We consider a linear approximation $F(\varphi) \equiv g(z)\varphi$ by freezing the function $g$ at some $z\in L^2(Q)$. We obtain the following linearised system by replacing the nonlinear term by above linear approximation,
\begin{align}
    \psi' + (\Delta)^2\psi - \Delta(g(z)\psi) &= U\chi_q, \text{ in } \Omega \times(0,T),\label{eq93a}\\
    \frac{\partial\psi}{\partial\textbf{n}}= \frac{\partial\Delta\psi}{\partial\textbf{n}} &= 0, \text{ on } \partial\Omega \times(0,T),\label{eq93b}\\ 
    \psi(0) &= \psi_0, \text{ in } \Omega.\label{eq93c}
\end{align}
Observe that if we replace $z$ by $\psi$, the above system coincides with the original system. 
 Note that, in the following analysis, we denote the term $g(z(x,t))$ by $a(x,t)$ for convenience.

First, we study the following control problem for the linearised system; for fixed time $T >0$ and for a given $\epsilon > 0$, $\omega = (\omega_1,..., \omega_k) \in \mathbb{R}^k$, whether there exists a control $U\in L^2(q)$ that steers the solution of \eqref{eq93a}-\eqref{eq93c} into $\epsilon$ neighbourhood of $\omega$ at points $x_i \in \mathcal{D}$ and time $T$, i.e., 
\begin{equation}\label{eq92}
    \sum_{i=1}^{k}(\psi(x_i,T)-\omega_i)^2 \leq \epsilon. 
\end{equation}

Now we state the existence of a strong solution of the linearised system, so that the pointwise evaluation of the solution in \eqref{eq92} makes sense.

\begin{theorem}[Existence of solution of the linearised system]\label{thm1.1}
    Let $\psi_0 \in H^2(\Omega)$, $a\in L^\infty(0,T; W^{1, \infty}(\Omega) \cap H^2(\Omega))$ and $U \in L^2(q)$. Then there exists a solution of the linearised system \eqref{eq93a}-\eqref{eq93c} with  regularity,
    \begin{align*}
        \psi \in \mathcal{V}:= L^\infty(0,T;H^2) \cap L^2(0,T; H^4)\cap H^1(0,T;H).
    \end{align*}
 Moreover, the solution satisfies the following energy estimate.
    \begin{equation}\label{eq2.9}
        \|\psi\|_{L^\infty(0,T;H^2(\Omega)) \cap L^2(0,T;H^4(\Omega)) \cap H^1(0,T;H)}^2 \leq C( \|\psi_0\|_{H^2(\Omega)}^2 + \|U\|_{L^2(q)}^2).
    \end{equation}
\end{theorem}

\begin{proof}
    The proof is by employing Galerkin approximation method; by deriving a priori estimates of the solution. Since the proof is similar to existence results proved in the literature \cite[Theorem 1, Theorem 2]{LKW}, we skip the proof. Moreover, the H\"older continuity of $\psi$ in the space variable is obtained by an application of  Aubin-Lions lemma, $L^\infty(0,T;H^2)\cap H^1(0,T;H) \hookrightarrow \hookrightarrow C([0,T];C^{0,\gamma}(\Bar{\Omega}))$, for $\gamma \in (0,1)$.  
\end{proof}

The above regularity of the linearised variable ensures the well-definedness of the problem described via \eqref{eq92}. We study the approximate controllability problem by a variational approach. In the following subsection, we introduce an adjoint system, and a cost functional in the adjoint variable. We reformulate the approximate controllability problem of the linearised system as a minimisation problem in the adjoint. Finally, we deduce a control which solves our approximate control problem in terms of the solution of the adjoint system associated with the minimiser of the cost. 

\subsubsection{The adjoint problem}
We define a cost functional, $\mathcal{J}: \mathbb{R}^k \longrightarrow \mathbb{R}$ given as follows, for $f =(f_1, f_2,..., f_k) \in \mathbb{R}^k$,
\begin{equation}\label{eq94}
    \mathcal{J}(f) = \mathcal{J}(f, g(z), \omega, \epsilon) = \underbrace{\frac{1}{2}(\int_q |\xi(x,t)|dxdt)^2}_{j_1} + \underbrace{\epsilon\big(\sum_{i=1}^{k}(f_i)^2\big)^{1/2} - \sum_{i=1}^{k}\omega_i f_i}_{j_2},
\end{equation}
where $\epsilon>0$, $\omega \in \mathbb{R}^k$ are fixed as in the approximate control problem and $\xi$ denotes the solution of the adjoint system,
\begin{align}
    -\xi' + (\Delta)^2\xi - a(x,t)\Delta\xi &= 0, \text{ in } \Omega\times(0,T),\label{eq95a}\\
    \frac{\partial\xi}{\partial\textbf{n}}= \frac{\partial\Delta\xi}{\partial\textbf{n}} &= 0,\text{ on } \partial\Omega\times(0,T),\label{eq95b}\\ 
    \xi(.,T) &= \xi_T, \text{ in } \Omega,\label{eq95c}
\end{align}
with terminal data, $\xi_T = \sum_{i=1}^{k}\delta_{x_i}f_i$. Recall that we are denoting $g(z(x,t))$ by $a(x,t)$. To motivate the above choice of cost functional, observe that if we establish the existence of a unique minimiser $\Hat{f}$ of $\mathcal{J}$,  by comparing the optimality condition, namely; $0 \in \partial \mathcal{J}(\Hat{f})$ and the variational equation satisfied by the linearised variable, we get a control which will imply the condition \eqref{eq92} and hence approximate controllability (as described above) of the linearised system.

Now we study the existence of a solution of the above adjoint problem. The solution is defined by transposition due to the low regular terminal data. 

\begin{definition}[Weak solution of the adjoint system by transposition]\label{def1}
     We call $\xi \in L^2(Q)$ a transposition solution of the adjoint system \eqref{eq95a}-\eqref{eq95c}, if  the following variational equation is satisfied.
  \begin{equation}\label{eq97}
        \int_{q}\xi U dxdt = {}_{\mathcal{M}(\Omega)}\langle \xi_T, \psi(.,T) \rangle_{C(\Bar{\Omega})}.
    \end{equation} 
    $\forall \; U \in L^2(q)$; which acts as a forcing term in the dual system (of the adjoint) given by,
    \begin{align}
    \psi' + (\Delta)^2\psi - \Delta(a(x,t)\psi) &= U\chi_{q}, \text{ in } \Omega \times(0,T),\label{eq96a}\\
    \frac{\partial\psi}{\partial\textbf{n}} = \frac{\partial\Delta\psi}{\partial\textbf{n}} &= 0, \text{ on } \partial\Omega \times(0,T),\label{eq96b}\\
    \psi(0) &= 0, \text{ in } \Omega.\label{eq96c}
\end{align}
\end{definition}

\begin{theorem}[Existence of a transposition solution of the adjoint system]\label{theorem1.2}
    Let $a\in L^\infty(0,T; W^{1, \infty}(\Omega) \cap H^2(\Omega))$, Then the following results hold. 
    \begin{enumerate}[label=\roman*)]
   \item There exists a transposition solution of the system \eqref{eq95a}-\eqref{eq95c} with $\xi_T = \sum_{i=1}^{k}\delta_{x_i}f_i$ in the sense of definition \ref{def1}. 
    \item The solution additionally satisfies,
    $$\xi \in L^{r'}(0,T; W^{2,p'}(\Omega)) \cap W^{1,r'}(0,T; W^{-2, p'}(\Omega)),$$
    where $U$ is considered to be in $L^p(0,T; W^{-2,p}(\Omega))$ for $1\leq r< \infty$, $2< p< \infty$ satisfying the condition, $2-\frac{4}{r}-\frac{2}{p} > 0$. Additionally, the following energy inequality holds.
    \begin{align}\label{eq2.15}
        \|\xi\|_{L^{r'}(0,T; W^{2,p'}(\Omega)) \cap W^{1,r'}(0,T; W^{-2, p'}(\Omega))} \leq C\|\xi_T\|_{\mathcal{M}(\Omega)},
    \end{align}
    where $C>0$ is a constant that depends on $a$ and $\Omega$.
    \end{enumerate}
    \end{theorem}

\begin{proof}
    proof of i) The proof is based on an approximation technique. Consider approximations to the terminal data, $\xi_T^n \in C(\Bar{\Omega})$ such that $\xi_T^n \overset{\ast}{\rightharpoonup} \xi_T$ in $\mathcal{M}(\Omega)$ and satisfies, 
   \begin{equation*}
       \|\xi_T^n\|_{L^1(\Omega)} \leq \|\xi_T\|_{\mathcal{M}(\Omega)}.
   \end{equation*}
 Let $\xi^n$ denote the solution of the following approximated system.  
\begin{align}
    -(\xi^n)' + (\Delta)^2\xi^n - a(x,t)\Delta\xi^n &= 0, \text{ in } \Omega\times(0,T),\label{eq2.10}\\
    \frac{\partial\xi^n}{\partial\textbf{n}}= \frac{\partial\Delta\xi^n}{\partial\textbf{n}} &= 0,\text{ on } \partial\Omega\times(0,T),\\ 
    \xi^n(.,T) &= \xi^n_T, \text{ in } \Omega,
\end{align}
By testing \eqref{eq2.10} with $U \in L^2(Q)$ we get,
   \begin{align}\label{eq2.018}
        \int_{Q} U\xi^n dxdt = \int_\Omega \psi(x,T)\xi_T^n dx.
   \end{align}
  where $\psi$ is the solution of system \eqref{eq96a}-\eqref{eq96c} with source term $U \in L^2(Q)$. Note that by theorem \ref{thm1.1} for $U\in L^2(Q)$, $\psi$ satisfies an energy inequality \eqref{eq2.9}.   
We apply H\"older's inequality to estimate the R.H.S of \eqref{eq2.018} and by \eqref{eq2.9}, we get, 
\begin{align*}
    |\int_{Q} U\xi^n dxdt| &\leq \|\psi(.,T)\|_{L^\infty(\Omega)}\|\xi_T^n\|_{L^1(\Omega)} 
    \leq \|\psi\|_{C([0,T];H^2(\Omega))}\|\xi_T\|_{\mathcal{M}(\Omega)} \leq C\|U\|_{L^2(q)}\|\xi_T\|_{\mathcal{M}(\Omega)}.
\end{align*}
Since $U$ is arbitrary, solution of the approximated system satisfies,
\begin{equation}\label{eq2.18}
    \|\xi^n\|_{L^2(Q)} \leq C\|\xi_T\|_{\mathcal{M}(\Omega)}.
\end{equation}
the R.H.S term in the above inequality is uniform in $n$, therefore by applying compactness arguments we can extract a subsequence,
\begin{align*}
    \xi^n \rightharpoonup \xi \text { in } L^2(Q).
\end{align*}
This enables us to recover \eqref{eq97} by passing to the limit as $n\rightarrow \infty$ in \eqref{eq2.018} with $\text{supp}(U) \subseteq q$, hence a transposition solution. 
\\Proof of ii)
Furthermore, we derive a higher regularity of the solution, namely, $\xi \in L^{r'}(0,T; W^{2,p'}(\Omega))$. To this end we consider $U \in L^r(0,T; W^{-2,p}(\Omega))$.  Note that, the solution of the linear system \eqref{eq96a}-\eqref{eq96c} with $U \in L^r(0,T; W^{-2,p}(\Omega))$ satisfies,
$$\psi \in L^r(0,T; W_N^{2,p}(\Omega)) \cap W^{1,r}(0,T; W_N^{-2,p}  (\Omega)) \hookrightarrow C([0,T]; (W^{-2,p}(\Omega), W^{2,p}(\Omega))_{1-\frac{1}{r}, r}).$$
\\This regularity will be discussed in detail in section 3, remark \ref{remark3.3}. In addition, we assume $1\leq r< \infty$, $2< p< \infty$ satisfies $2-\frac{4}{r}-\frac{2}{p} > 0$, then  we have $\psi \in C([0,T]; C_b(\Omega))$ (See remark \ref{remark2}).

We say $\xi \in L^{r'}(0,T; W^{2,p'})$ is a solution of the adjoint system if $ \forall U \in L^r(0,T;W^{-2,p}(\Omega))$, the following identity holds.
\begin{equation}\label{eq2.28}
        \int\limits_{0}^{T}{}_{W^{2,p'}(\Omega)}\langle\xi, U \rangle_{W^{-2,p}(\Omega)} = {}_{\mathcal{M}(\Omega)}\langle \xi_T, \psi(.,T) \rangle_{C(\Bar{\Omega})}.
    \end{equation} 
Consider analogous identity satisfied by the approximated solution say $\xi^n$ associated with the terminal data $\xi_T^n$. Applying H\"olders to estimate the R.H.S of the above inequality leads to, 
\begin{align*}
    |\int\limits_{0}^{T}\langle\xi^n, U \rangle | &\leq \|\xi_T^n\|_{L^1(\Omega)}\|\psi(.,T)\|_{L^\infty(\Omega)} \leq \|\xi_T\|_{\mathcal{M}(\Omega)}\|\psi(.,T)\|_{L^\infty(\Omega)} \\
    &\leq \|\xi_T\|_{\mathcal{M}(\Omega)}\|U\|_{L^r(0,T;W^{-2,p}(\Omega))}
\end{align*}
The last inequality follows from \eqref{eq3.15} and the embedding given in remark \ref{remark2}. Thus,
\begin{align*}
    \|\xi^n\|_{L^{r'}(0,T; W^{2,p'}(\Omega))} &\leq \|\xi_T\|_{\mathcal{M}(\Omega)}.
\end{align*}
Therefore by compactness arguments we obtain a weak limit $\xi\in L^{r'}(0,T; W^{2,p'}(\Omega))$ which satisfies \eqref{eq2.28}. Observe that, $\xi\in L^{r'}(0,T; W^{2,p'}(\Omega)) \hookrightarrow L^{r'}(0,T; C(\Bar{\Omega}))$ by Sobolev embedding. 

Further, the regularity in time, $\xi \in L^{r'}(0,T; W^{-2, p'}(\Omega))$ can be derived by testing \eqref{eq95a} with an arbitrary test function $\rho \in L^r(0,T; W^{2,p}(\Omega))$. Consider, 
\begin{align*}
    -\langle \xi', \rho \rangle + (\Delta\xi, \Delta \rho) - (a(x,t)\Delta\xi, \rho) &= 0.
\end{align*}
By applying H\"olders and Young's inequalities we obtain,
\begin{align*}
    |\langle \xi', \rho \rangle| &\leq \|\Delta\xi\|_{L^{p'}}\|\Delta \rho\|_{L^p} + \|a\|_{L^\infty}\|\Delta\xi\|_{L^{p'}}\|\rho\|_{L^p}\\
    &\leq (\|\Delta\xi\|_{L^{p'}} + \|a\|_{L^\infty}\|\Delta\xi\|_{L^{p'}}\|)\|\rho\|_{W^{2,p}}.
\end{align*}
This yields the regularity in time, $\xi \in L^{r'}(0,T;W^{-2,p'}(\Omega))$. 
\end{proof}

\begin{remark}
        Note that, to study the approximate controllability problem for the semilinear equation with a distributional $L^2$ control, we will only require the result \romannumeral 1) in the theorem \ref{theorem1.2}. The additional regularity of $\xi$ in \romannumeral 2) of theorem \ref{theorem1.2} is used to study the approximate controllability problem associated with the pointwise control in section 3.  Hence the detailed discussion in the proof is deferred to section 3.
 \end{remark}

\begin{remark} 
Consider the map $f \mapsto \xi_f$ from $\mathbb{R}^k \rightarrow L^1(q)$, where $\xi_f$ denotes the solution of the adjoint system corresponding to the terminal data, $\xi(.,T) = \sum_{i=1}^{k}\delta_{x_i}f_i$. For each $f\in\mathbb{R}^k$, the adjoint system admits a solution $\xi_f\in L^2(q)\subset L^1(q)$. Hence, the above map is well defined, and consequently, the cost functional $\mathcal{J}$ in \eqref{eq94} is well defined. Moreover, this map is linear and continuous.
\end{remark}

The above regularity of the adjoint is not sufficient to study the minimisation problem of the cost $\mathcal{J}$. In the following we state and prove a regularity theorem  for the adjoint with a low regular terminal term. The result is achieved by deriving higher regularity of the approximated system. However, the higher regularity is obtained only on a time interval away from the terminal time. The following regularity is necessary to deduce the unique continuation property of the adjoint which is crucial in the minimisation problem.

\begin{theorem}\label{thm2.3}
    The solution of the adjoint problem \eqref{eq95a}-\eqref{eq95c} has the following additional regularity for each $0< \delta < T$.
    \begin{align*}
        \xi \in L^2(0,T-\delta, H^4(\Omega))\cap L^\infty(0,T-\delta; H^2(\Omega)) \cap  H^1(0,T-\delta;H^{-2}(\Omega)).
    \end{align*}
\end{theorem}
\begin{proof}
Consider the approximations to the terminal data $\xi_T \in \mathcal{M}(\Omega)$ by functions $\xi_T^n \in C(\Bar{\Omega})$ introduced in the previous theorem.
We derive higher regularity of $\xi$, the solution corresponding to $\xi_T$, by deriving a priori estimates for $\xi^n$ uniform in $n$. Consider the following equation obtained by testing \eqref{eq2.10} with $(T-t)\xi^n$,  
\begin{align}\label{eq2.20}
   -\langle(\xi^n)', (T-t)\xi^n\rangle + (T-t)\|\Delta\xi^n\|^2 = -(T-t)(a\Delta\xi^n, \xi^n).
\end{align}
Observe that 
\begin{align*}
    -\langle(\xi^n)', (T-t)\xi^n\rangle &= -\frac{1}{2}\frac{d}{dt}\big((T-t)\|\xi^n\|^2\big) - \frac{1}{2}\|\xi^n(t)\|^2, \\
    |(T-t)(a\Delta\xi^n, \xi^n)| &\leq \frac{(T-t)}{2}\|\Delta\xi^n\|^2 + \frac{\|a\|_{L^\infty}(T-t)}{2}\|\xi^n\|^2.
\end{align*}
Substituting back in \eqref{eq2.20}, we get,
\begin{align}
    -\frac{1}{2}\frac{d}{dt}\big((T-t)\|\xi^n\|^2\big) + \frac{(T-t)}{2}\|\Delta\xi^n\|^2 &\leq \frac{(1+\|a\|_{L^\infty}T)}{2}\|\xi^n\|^2. \label{eq2.19} 
\end{align}
Integrating \eqref{eq2.19} over the interval $[t, T]$ results in,
\begin{align}\label{eq2.024}
    \frac{1}{2}(T-t)\|\xi^n(t)\|^2 + \int\limits_{t}^{T} \frac{(T-\tau)}{2}\|\Delta\xi^n(\tau)\|^2 d\tau &\leq C_T\int\limits_{t}^{T}\|\xi^n(s)\|^2ds, 
\end{align}
where $C_T = \frac{(1+\|a\|_{L^\infty(Q)}T)}{2}$.
Observe that the term on R.H.S of the last inequality is bounded uniformly in $n$ using the estimate \eqref{eq2.18}. Note that, $(T-t)$ is bounded below by $\delta$ when $t\in [0,T-\delta]$. Consequently,
\begin{align}\label{eq2.21}
    \|\xi^n(t)\|^2 &\leq \frac{2}{\delta}C(1+T)\|\xi_T\|_{\mathcal{M}(\Omega)}^2, \text{ for } 0\leq t\leq T-\delta,
\end{align}
where $C >0$ is a constant. Further by integrating \eqref{eq2.19} over $[0, T-\delta]$ and applying the estimate \eqref{eq2.18}, \eqref{eq2.21}, we arrive at,
\begin{align}
  \frac{T}{2}\|\xi^n(0)\|^2 + \int\limits_{0}^{T-\delta}\frac{(T-\tau)}{2}\|\Delta\xi^n(\tau)\|^2 d\tau &\leq \frac{\delta}{2}\|\xi^n(T-\delta)\|^2 + C\int\limits_{0}^{T-\delta}\|\xi^n(\tau)\|^2 d\tau \nonumber\\ 
  &\leq  \frac{\delta}{2}\|\xi^n(T-\delta)\|^2 +  C \|\xi_T\|_{\mathcal{M}(\Omega)}^2 \nonumber\\
  &\leq  C(1+T) \|\xi_T\|_{\mathcal{M}(\Omega)}^2. \nonumber\\
 \text{Hence, \;\; } \int\limits_{0}^{T-\delta}\|\Delta\xi^n(\tau)\|^2 d\tau &\leq \frac{2}{\delta}C(1+ T)\|\xi_T\|_{\mathcal{M}(\Omega)}^2. \label{eq2.22}
\end{align}
By exploiting the elliptic regularity of $\xi^n$ and uniform estimates \eqref{eq2.21}, \eqref{eq2.22} we get, $\xi^n$ is uniformly bounded in $L^\infty(0,T-\delta; H) \cap L^2(0,T-\delta; H^2(\Omega))$ for fixed $\delta \in (0,T)$. By applying compactness arguments, we obtain, for any $\delta \in (0,T)$ 
$$ \xi \in L^\infty(0,T-\delta; H) \cap L^2(0,T-\delta; H^2(\Omega)).$$
Further by testing \eqref{eq95a} with an arbitrary $\rho \in L^2(0,T-\delta;H^2(\Omega))$, we obtain,
\begin{align*}
    \langle \xi', \rho \rangle + (\Delta\xi, \Delta\rho) + (\nabla(a(x,t)\xi), \nabla\rho) = \langle U, \rho \rangle.
\end{align*}
Estimating each of its terms by applying H\"olders results,
\begin{align*}
    |\langle \xi', \rho \rangle|  &\leq \|\Delta\xi\|\|\Delta\rho\| + \|\nabla a\|_{L^\infty}\|\|\xi\|\|\Delta\rho\|\\
    &\leq \big(\|\Delta\xi\| + \|\nabla a\|_{L^\infty}\|\|\xi\| \big)\|\rho\|_{H^2}
\end{align*}
Integrating both sides over an interval $(0, T-\delta)$ and applying estimates \eqref{eq2.22}, \eqref{eq2.18} yields the following regularity. 
$$\xi' \in L^2(0,T-\delta; H^{-2}(\Omega)).$$

To derive a higher regularity of $\xi^n$ we test \eqref{eq2.10} with $(T-t)^2(\Delta)^2 \xi^n$. This yields,
\begin{align}\label{eq2.022}
   -\langle(\xi^n)', (T-t)^2(\Delta)^2\xi^n\rangle + (T-t)^2\|(\Delta)^2\xi^n\|^2 = -((T-t)a\Delta\xi^n, (T-t)(\Delta)^2\xi^n).
\end{align}
Observe that 
\begin{align*}
    -\langle(\xi^n)', (T-t)^2(\Delta)^2\xi^n\rangle &= -\frac{1}{2}\frac{d}{dt}\big((T-t)^2\|\Delta\xi^n\|^2\big) - \frac{1}{2}(T-t)\|\Delta\xi^n(t)\|^2, \\
    |((T-t)a\Delta\xi^n, (T-t)(\Delta)^2\xi^n)| &\leq \frac{(T-t)^2}{2}\|(\Delta)^2\xi^n\|^2 + \frac{\|a\|_{L^\infty}(T-t)^2}{2}\|\Delta\xi^n\|^2.
\end{align*}
Substituting back in \eqref{eq2.022}, we get,
\begin{align}\label{eq2.23}
    -\frac{1}{2}\frac{d}{dt}\big((T-t)^2\|\Delta\xi^n\|^2\big) + \frac{(T-t)^2}{2}\|(\Delta)^2\xi^n\|^2
    &\leq \frac{(T-t)}{2}\|\Delta\xi^n(t)\|^2 + \frac{(T-t)^2\|a\|_{L^\infty}}{2}\|\Delta\xi^n\|^2.
\end{align}
Integrate the above inequality over $[t,T]$, we obtain,
\begin{align*}
    \frac{1}{2}(T-t)^2\|\Delta\xi^n(t)\|^2  &\leq  \int\limits_{t}^{T}\frac{(T-\tau)}{2}\|\Delta\xi^n(\tau)\|^2 d\tau + \|a\|_{L^\infty}\int\limits_{t}^{T}\frac{(T-\tau)^2}{2}\|\Delta\xi^n(\tau)\|^2 d\tau.
\end{align*}
Now observe from \eqref{eq2.024}, \eqref{eq2.18} that the first term on the R.H.S of the above inequality is uniformly bounded in n. Thus, 
\begin{align*}
    \frac{1}{2}(T-t)^2\|\Delta\xi^n(t)\|^2  &\leq  C(1+T)\|\xi_T\|_{\mathcal{M}(\Omega)}^2 + \|a\|_{L^\infty}\int\limits_{t}^{T}\frac{(T-\tau)^2}{2}\|\Delta\xi^n(\tau)\|^2 d\tau.
\end{align*}
By exploiting Gronwall's inequality we obtain,
\begin{align*}
    (T-t)^2\|\Delta\xi^n(t)\|^2 &\leq C(1+T)\|\xi_T\|_{\mathcal{M}(\Omega)}^2 \text{exp}\Big( -\int\limits_{0}^{t}\|a\|_{L^\infty}ds\Big)\\
    &\leq C(1+T)\|\xi_T\|_{\mathcal{M}(\Omega)}^2 e^{-\|a\|_{L^\infty}t}.
\end{align*}
If $t\in [0,T-\delta]$, we have $(T-t) \geq \delta$, therefore,
\begin{equation}\label{eq2.26}
    \|\Delta\xi^n(t)\|^2 \leq \frac{1}{\delta^2}C(1+T)\|\xi_T\|_{\mathcal{M}(\Omega)}^2, \text{ for } 0\leq t \leq T-\delta.
\end{equation}
Moreover, by integrating \eqref{eq2.23} over $[0,T-\delta]$, and applying estimates \eqref{eq2.22}, \eqref{eq2.26}, we get,
\begin{align}
    \frac{T^2}{2}\|\Delta\xi^n(0)\|^2 + \int\limits_{0}^{T-\delta}\frac{(T-\tau)^2}{2}\|(\Delta)^2\xi^n(\tau)\|^2 d\tau
    &\leq \frac{\delta^2}{2}\|\Delta\xi^n(T-\delta)\|^2  + \frac{T}{2}\int\limits_{0}^{T-\delta}\|\Delta\xi^n(t)\|^2 + \frac{T^2\|a\|_{L^\infty}}{2}\int\limits_{0}^{T-\delta}\|\Delta\xi^n\|^2 \nonumber\\
    &\leq C(1+T)\|\xi_T\|_{\mathcal{M}(\Omega)}^2 + CT(1+T)\int\limits_{0}^{T-\delta}\|\Delta\xi^n\|^2 \nonumber \\
    &\leq C(1+C_T)(1+\frac{1}{\delta})\|\xi_T\|_{\mathcal{M}(\Omega)}^2, \nonumber\\
    \implies \int\limits_{0}^{T-\delta}\|(\Delta)^2\xi^n(\tau)\|^2 d\tau &\leq \frac{1}{\delta^2}(1+\frac{1}{\delta})C(1+C_T)\|\xi_T\|_{\mathcal{M}(\Omega)}^2. \label{eq2.27}
\end{align}

By exploiting the elliptic regularity of $\xi^n$, using estimates \eqref{eq2.22}, \eqref{eq2.27} we obtain, $\xi^n$ is uniformly bounded in $L^2(0,T-\delta; H^4(\Omega)) \cap L^\infty(0,T-\delta; H^2(\Omega))$. Moreover, by compactness arguments we can extract a convergent subsequence of $\xi^n$. it follows that the limit,
$$ \xi \in L^2(0,T-\delta; H^4(\Omega)) \cap L^\infty(0,T-\delta; H^2(\Omega)) $$
solves the system \eqref{eq95a}-\eqref{eq95c}. This completes the proof. 
\end{proof}

Our aim is to show that $\mathcal{J}$ has a unique minimiser. To this end, we prove $\mathcal{J}$ is continuous, strictly convex and coercive. Using continuity of the map $f\mapsto \xi_f$, we can conclude that $\mathcal{J}$ is continuous. 
Further we prove that $\mathcal{J}$ is strictly convex, observe that $j_2$ as denoted in \eqref{eq94} is convex. Let $f, h \in \mathbb{R}^k$ and $\xi_f, \xi_h$ be respective solutions to the adjoint problem \eqref{eq95a} - \eqref{eq95c}, then $j_1$ satisfies,
\begin{align*}
    j_1(\theta f +(1-\theta)h) &= \frac{1}{2}\|\theta\xi_f+ (1-\theta)\xi_h\|_{L^1(q)}^2 \leq \frac{1}{2}\big(\theta\|\xi_f\|_{L^1(q)}+ (1-\theta)\|\xi_h\|_{L^1(q)}\big)^2 \\
   &\leq \frac{1}{2}\big( \theta\|\xi_f\|_{L^1(q)}^2 + (1-\theta)\|\xi_h\|_{L^1(q)}^2\big) = \theta j_1(f) + (1-\theta) j_1(h).
\end{align*}
The second last inequality follows by the convexity of the function $x^2$ on $[0, \infty)$. Therefore, $j_1$ is convex.

Additionally, we prove the cost $\mathcal{J}$ is strictly convex. Since the last term of $\mathcal{J}$ in \eqref{eq94} is linear in $f$, we get the strict convexity condition from the first two terms. Observe that, even if one of the terms satisfies a strict convexity condition then the cost $\mathcal{J}$ is strictly convex. For the second term, the equality condition of the H\"older's inequality holds whenever $h = \alpha f$, for a constant $\alpha$. Hence this term evaluated at $ (\theta f + (1- \theta) h)$  will be strictly convex  unless $h = \alpha f$. Now if $h = \alpha f$ we need to show that the strict convexity holds for $j_1$. This follows from the observation that the equality in convexity condition of $j_1$ implies $\xi_{f}\xi_{h} \geq 0$ and $\|\xi_{f}\|_{L^1(q)} = \|\xi_{h}\|_{L^1(q)}$. Therefore we conclude that $\mathcal{J}$ is strictly convex.  

Now, we prove $\mathcal{J}$ is coercive. i.e. there exists a constant $\alpha >0 $ such that,
\begin{align}\label{eq2.17}
    \liminf\limits_{\|f\|\rightarrow \infty} \frac{\mathcal{J}(f, a, \omega)}{\|f\|} \geq \alpha.
\end{align}
We use the unique continuation property of the adjoint system to prove the coercivity of $\mathcal{J}$. 

\begin{theorem}[Unique continuation property of the adjoint system]\label{theorem2.4}
    Let $q$ be a non-empty open subset of $Q$. Suppose the solution of the adjoint system, $\xi_f \equiv 0$ on $q$, then $f_i = 0$ for $1\leq i\leq k$.
\end{theorem}
\begin{proof}
We have from theorem \ref{thm2.3}, $\xi \in L^2(0,T-\delta; H^4(\Omega))$, for any $\delta \in (0,T)$. Then the unique continuation property of $\xi$ in $ \Omega \times (0, T-\delta)$ follows immediately from \cite[Theorem 3.2]{saut87}. If $\xi = 0$ a.e in $q \subseteq Q$, nonempty, open;  then $\xi = 0$ a.e in $\Omega \times (0, T-\delta)$ for a $\delta \in (0,T)$. 
\\Since $\delta$ is arbitrary we get, $\xi = 0$ a.e in $Q$. Therefore, we have from the definition  \ref{def1}, of a solution,  for any $U\chi_{q} \in L^2(q)$,
\begin{align}\label{eq2.026}
 {}_{\mathcal{M}(\Omega)}\langle \xi_T, \psi(.,T) \rangle_{C(\Bar{\Omega})} = \int_{q}\xi U dxdt = 0.
\end{align}
\\Let $\rho \in \mathcal{D}(\Omega)$, We define $\psi(x,t) = \frac{t}{T}\rho$. Then $\psi$ satisfies \eqref{eq96a}-\eqref{eq96c} for $U = \frac{1}{T}\rho + (\Delta)^2\rho + \Delta(a \rho)$. Indeed, $U \in L^2(Q)$ because $a \in L^\infty(W^{1,\infty}(\Omega) \cap H^2(\Omega))$. Therefore the image set of the solution map $U \mapsto \psi(.,T)$ from $L^2(Q) \rightarrow H^2(\Omega)$ defined by theorem \ref{thm1.1} is dense in $C(\Bar{\Omega})$. 
Since \eqref{eq2.026} holds  for every $\psi (., T) =\rho \in \mathcal{D}(\Omega)$,
it follows that, $\xi_T = 0$ in $\mathcal{M}(\Omega)$. 
This implies that whenever the transposition solution of the system \eqref{eq95a}- \eqref{eq95c} corresponding to $\xi_T = \sum_{i=1}^{k}f_i\delta_{x_i}$ denoted by  $\xi$ is equal to 0 in $q$ then $f_i = 0$ for each $i$, $1\leq i \leq k$.  Hence unique continuation property holds.
\end{proof}
\noindent Consequently, we obtain the following result on the coercivity of $\mathcal{J}$.

\begin{coro}
    The cost functional $\mathcal{J}$, defined in \eqref{eq94} is coercive. In particular, 
  \begin{align}\label{eq93.22}
    \liminf\limits_{\|f\|\rightarrow \infty} \frac{\mathcal{J}(f, a, \omega)}{\|f\|} \geq \epsilon, 
    \end{align}
    where $\epsilon$ is the constant in the definition of the cost $\mathcal{J}$ given in \eqref{eq94}.
\end{coro}
\begin{proof}
    The proof is by contradiction. Suppose to the contrary, there exists a sequence $(f^n)_{n\geq 1}$ in $\mathbb{R}^k$ such that $\|f^n\| \rightarrow \infty$ and $\lim\inf_{n\rightarrow \infty}\frac{\mathcal{J}(f^n)}{\|f^n\|} <\epsilon$. 
    We introduce the normalised variable $\Tilde{f}^n = \frac{f^n}{\|f^n\|}$. Since $\Tilde{f}^n$ are uniformly bounded we can extract a subsequence still denoting by same index, $\Tilde{f}^n$ such that, $\Tilde{f}^n \rightarrow \Tilde{f}$ in $\mathbb{R}^k$. With the notation solution corresponding to $f$ is denoted by $\xi_{f}$, observe that,
\begin{align*}
    \frac{\mathcal{J}(f^n)}{\|f^n\|} &= \frac{\|f^n\|}{2}(\int_q |\frac{\xi_{f^n}}{\|f^n\|}|)^2 + \epsilon - \sum_{i=1}^{k}\omega_i\Tilde{f}^n_i.
\end{align*}
Since $\frac{\xi_{f^n}}{\|f^n\|} = \xi_{\frac{f^n}{\|f^n\|}} = \xi_{\Tilde{f}^n}$. Therefore,
\begin{align*}
    \frac{\mathcal{J}(f^n)}{\|f^n\|} &= \frac{\|f^n\|}{2}(\int_q |\xi_{\Tilde{f}^n}|)^2 + \epsilon - \sum_{i=1}^{k}\omega_i\Tilde{f}^n_i. 
\end{align*}
Observe that, $\epsilon - \sum_{i=1}^{k}\omega_i\Tilde{f}^n_i $ is finite as $\epsilon$, $\omega$ are fixed and $\Tilde{f}^n$ is a unit vector. Since $\|f^n\| \rightarrow \infty$ and by our assumption on the limit of $\mathcal{J}(f^n)/\|f^n\|$, we get the limit $\lim_{n\rightarrow \infty}(\int_q |\xi_{\Tilde{f}^n}|) =0$. Therefore, $\xi_{\tilde{f}} \equiv 0$ in $q$ and by unique continuation property $\Tilde{f} =0$ in $\mathcal{M}(\Omega)$. Hence a contradiction to the fact that $\Tilde{f}$ is a limit of unit vectors. This completes the proof.
\end{proof}

Since $\mathcal{J}$ is continuous, strictly convex and coercive, it follows that $\mathcal{J}$ has a unique minimiser in $\mathbb{R}^k$. Let $\Hat{f} \in \mathbb{R}^k$ denotes the minimiser of $\mathcal{J}(f, g(z), \epsilon, \omega)$. Recall the notation $g(z) = a(x,t)$ for $z\in L^2(Q)$, which will be again used in subsection 2.2 to study the nonlinear problem. We define a map, \begin{align*}
    \mathcal{N}: L^\infty(Q) \times \mathbb{R}^k &\rightarrow \mathbb{R}^k \text{ such that}\\
   (g(z),\omega) &\mapsto \Hat{f}.
\end{align*}
where $\hat f$ is the minimiser of $\mathcal{J}$, for fixed $\epsilon >0$ and $\omega \in \mathbb{R}^k$. We can prove that if $B\subseteq L^\infty(Q)$ is bounded and $K\subseteq \mathbb{R}^k$ is compact, then $\mathcal{N}(B\times K)$ is bounded subset of $\mathbb{R}^k$. The proof is given by establishing the constant $\epsilon$ in coercivity condition \eqref{eq93.22} is uniform in $g$ and $\omega$ \cite{Fabre95, DRA}. Indeed, this constant is same as the constant $\epsilon$ in the cost defined in \eqref{eq94}. We prove a lemma which will be used in the analysis of the nonlinear problem.

\begin{lemma}\label{prop2.5}
    Suppose, $\{z_n\}\subseteq L^2(Q), \{ \omega_n\} \subseteq \mathbb{R}^k$, $n\in \mathbb{N}$ be sequences such that,
    \begin{align}
        g(z_n) &\overset{\ast}{\rightharpoonup} g(z_0) \text{ in } L^\infty(Q),\label{eq1.23}\\
        \omega_n &\rightarrow \omega \text{ in } \mathbb{R}^k,\label{eq1.023}
    \end{align}
     Then we have, upto a subsequence in $\mathbb{R}^k,$ (denoted by same suffixes), $\mathcal{N}(g(z_n), \omega_n) \rightarrow \mathcal{N}(g(z_0), \omega)$. 
\end{lemma}

\begin{proof}
 We use the notation,
 \begin{align*}
     \hat{f}^n = \mathcal{N}(g(z_n), \omega_n), \hspace{.5cm} \hat{f} = \mathcal{N}(g(z_0), \omega).
 \end{align*}
Since $g(z_n)$, $\omega_n$ are uniformly bounded in $L^\infty(Q)$ and $ \mathbb{R}^k$ respectively, we have $\hat{f}^n$ uniformly bounded in $\mathbb{R}^k$ and $\exists f \in \mathbb{R}^k$ such that, upto a subsequence,
\begin{align*}
    \hat{f}^n \rightarrow f \text{ in } \mathbb{R}^k.
\end{align*}
Therefore by continuous dependence of solution on terminal data from \eqref{eq2.18}, we get,
\begin{align*}
    \xi_{\hat{f}^n} \rightarrow \xi_f \text{ in } L^1(q).
\end{align*}
Using weak lower semi continuity of $\mathcal{J}$, we obtain,
\begin{align*}
    \mathcal{J}(f, g(z_0), \omega) \leq \liminf_{n\rightarrow \infty}{\mathcal{J}(\hat{f}^n, g(z_n), \omega_n)}.
\end{align*}
Since $\hat{f}^n$ is the minimiser of $J(., g(z_n),\omega_n)$, observe that
\begin{align*}
    \mathcal{J}(\hat{f}^n, g(z_n), \omega_n) \leq \mathcal{J}(\hat{f}, g(z_n), \omega_n).
\end{align*}
Now by passing to the limit $n\rightarrow \infty$ on both sides of the above inequality, we get,
\begin{align*}
    \mathcal{J}(f, g(z_0), \omega_0) \leq \mathcal{J}(\hat{f}, g(z_0), \omega_0).
\end{align*}
Since $\hat{f} = \mathcal{N}(g(z_0), \omega)$, using the uniqueness of the minimiser, we get, $\hat{f} = f$. We can conclude that there exists a subsequence of $ \mathcal{N}(g(z_n), \omega_n)$  that converges to $\mathcal{N}(g(z_0), \omega)$ in $\mathbb{R}^k$. 
\end{proof}

Since the cost $\mathcal{J}$ is not differentiable, we use the notion of sub-differential to derive an optimality condition. In the following theorem, we characterise the sub-differential of $\mathcal{J}$ at a point $f \in \mathbb{R}^k$. We have at the minimiser, say $\Hat{f}\in \mathbb{R}^k$, $0\in \partial \mathcal{J}(\Hat{f})$, the sub-differential of $\mathcal{J}$ at $\Hat{f}$. This characterisation is useful to derive a control which solves our control problem. Note that, by $v \in \text{sgn}(h)$ we mean that $v(x,t) = \frac{h(x,t)}{|h(x,t)|}$ for $|h(x,t)|\neq 0$ and $|v(x,t)|\leq 1$ for $|h(x,t)| = 0$.

\begin{prop}\label{prop2.7}
    Let $f\in \mathbb{R}^k$, $f \neq 0$, The sub-differential of $\mathcal{J}$ at $f$ is given by,
    \begin{align}\label{partial J}
        \partial \mathcal{J}(f) = \bigg\{ \theta \in \mathbb{R}^k :\,\, &\exists v\in \text{sgn}(\xi_f) \text{ which satisfies }, \nonumber\\
        &(\theta, h) = \|\xi_f\|_{L^1(q)}\int_q v \xi_h + \epsilon\sum_{i=1}^{k}\frac{f_i h_i}{\|f\|} -\sum_{i=1}^{k}\omega_ih_i, \hspace{.5cm}\forall h \in \mathbb{R}^k \bigg\}.
    \end{align}
\end{prop}
\begin{proof}
Consider the cost $\mathcal{J}$ defined in \eqref{eq94}. Observe that $j_2$ is Fr\'etchet differentiable and for every $f \neq 0$,
\begin{equation*}
    ( D j_2(f), h) = \epsilon\sum_{i=1}^{k}\frac{f_i h_i}{\|f\|} -\sum_{i=1}^{k}\omega_i h_i, \hspace{.5cm} \forall h \in \mathbb{R}^k. 
\end{equation*}
To determine the sub-differential of $j_1$ at $f$, we have,  $\theta \in \partial j_1(f)$, iff $\forall h \in \mathbb{R}^k$, 
\begin{equation*}
    (\theta, h) \leq \lim_{s\rightarrow 0^+} \frac{j_1(f + sh)- j_1(f)}{s},
\end{equation*}
using $j_1$ is convex (See \cite{AUE}). And the limit,
\begin{align}
    \lim_{s\rightarrow 0^+} \frac{j_1(f + sh)- j_1(f)}{s} &= \lim_{s\rightarrow 0^+} \frac{1}{2}\bigg(\int_q |\xi_{f + sh}| + |\xi_f|dxdt \bigg)\bigg(\frac{\int_q |\xi_{f + sh}| - |\xi_f| dxdt }{s}\bigg) \nonumber\\
    &= \big(\int_q |\xi_f|dxdt\big)\lim_{s\rightarrow 0^+} \frac{\int_q |\xi_{f + sh}| - |\xi_f| dxdt}{s} \nonumber\\
    &= \|\xi_f\|_{L^1(q)}\bigg( \int_{q\cap \{(x,t) \in Q: \xi_f \neq 0\}}\text{sgn}(\xi_f)\xi_h dxdt + \int_{q\cap \{(x,t) \in Q: \xi_f = 0\}}|\xi_h| dxdt \bigg).\label{eq2.37}
\end{align}

Our aim is to characterise $\theta \in \partial j_1(f)$. To this end, we define a continuous linear functional $\Psi: \mathcal{R}\rightarrow \mathbb{R}$ by $\Psi(\xi_h):= (\theta, h)$. Here $\mathcal{R}$ represents the range of the solution map $h\mapsto \xi_h$. Since the map $\Psi$ is continuous linear, we can extend the map from $\mathcal{R}$ to $L^1(q)$ by Hahn-Banach theorem. Further by characterising the extended map, say $\Tilde{\Psi}$, we obtain an expression for the inner product $(\theta, h)$ which leads to a characterisation of $\partial j_1(f)$.

Consider the map $\Tilde{\Psi} : L^1(q) \rightarrow \mathbb{R}$ satisfying, 
\begin{enumerate}[label=\roman*)]
\item The restriction of $\Tilde{\Psi}$ to $\mathcal{R}$ agrees with $\Psi$.
     \begin{align*}
    \Tilde{\Psi}(\xi_h) &= (\theta, h), \hspace{.5cm}\forall h\in \mathbb{R}^k.
     \end{align*}
\item for every $\Phi \in L^1(q)$, the following inequality holds.
    \begin{align}\label{eq1.25}
    \Tilde{\Psi}(\Phi) &\leq \|\xi_f\|_{L^1(q)}\bigg( \int_{q\cap \{(x,t): \xi_f \neq 0\}}\text{sgn}(\xi_f)\Phi dxdt + \int_{q\cap \{(x,t): \xi_f = 0\}}|\Phi|dxdt \bigg).
   \end{align}
\end{enumerate}
It follows that,
\begin{align*}
    |\Tilde{\Psi}(\Phi)| \leq C\|\xi_f\|_{L^1(q)}\|\Phi\|_{L^1(q)}, \forall \Phi \in L^1(q),
\end{align*}
for some constant $C>0$. Therefore $\Tilde{\Psi}$ is a continuous linear functional on $L^1(q)$ which means $\Tilde{\Psi} \in L^1(q)'$. Therefore, $\Tilde{\Psi} \in L^\infty(q)$. Further, from \eqref{eq1.25}, it follows that, 
\begin{align*}
    \bigg|\int_q \Tilde{\Psi}(x,t)\Phi(x,t)dxdt - \|\xi_f\|_{L^1(q)}\int_{q\cap \{(x,t): \xi_f \neq 0\}}\text{sgn}(\xi_f)\Phi dxdt \bigg| &\leq \|\xi_f\|_{L^1(q)}\int_{q\cap \{(x,t): \xi_f = 0\}}|\Phi| dxdt.
\end{align*}

Consider the following two cases,
\begin{enumerate}
    \item If $\text{Supp}(\Phi) \subseteq q\cap \{(x,t)\in Q: \xi_f \neq 0\} $, Then we get $\Tilde{\Psi}(x,t) = \|\xi_f\|_{L^1(q)}\frac{\xi_f(x,t)}{|\xi_f(x,t)|}$. 
    \item If $\text{Supp}(\Phi) \subseteq q\cap \{(x,t)\in Q: \xi_f = 0\}$, Then we have $|\Tilde{\Psi}(x,t)| \leq \|\xi_f\|_{L^1(q)}$.
\end{enumerate}
From above two cases we conclude that $\Tilde{\Psi}(x,t) = \|\xi_f\|_{L^1(q)}v\chi_q$, $v\in \text{sgn}(\xi_f)$. Therefore, $\theta \in \partial j_1(f)$ implies $(\theta, h) = \|\xi_f\|_{L^1(q)}\int_q v \xi_h$, $\forall h\in \mathbb{R}^k$.
\\Conversely suppose, $\Tilde{\Psi}(x,t) = \|\xi_f\|_{L^1(q)}v\chi_q$ for some $v\in \text{sgn}(\xi_f)$, we define a map from $\mathbb{R}^k$ to $\mathbb{R}$ by $h \mapsto \int_q \Tilde{\Psi}(x,t)\xi_h$. This map is linear continuous. Therefore by Riesz representation theorem, there exists $\theta \in \mathbb{R}^k$ such that $(\theta, h) = \int_q \Tilde{\Psi}(x,t)\xi_h$. Moreover, comparing \eqref{eq2.37} and \eqref{eq1.25} we conclude that $\theta \in \partial j_1(f)$.
\end{proof}

The following theorem is on the pointwise approximate controllability of the linearised system. It is important to observe that if the desired state $\omega$ with $\|\omega\| \leq \epsilon$, then for the zero control, the solution $\varphi = 0$ satisfies the approximate controllability condition \eqref{eq92}. For the case $\|\omega\| > \epsilon$, we have the following result.

\begin{theorem}\label{thm1.7}
    Let $\|\omega\| > \epsilon$ and $\Hat{f}$ be the minimiser of $\mathcal{J}(., a, \omega)$. Then for $U= \|\xi_{\Hat{f}}\|_{L^1(q)} \emph{ sgn }\xi_{\Hat{f}} \chi_q$, the solution of linearised system \eqref{eq96a}-\eqref{eq96c} satisfies the approximate controllability condition \eqref{eq92}.
\end{theorem}
\begin{proof}
Consider $\Hat{f}$, the minimiser of $\mathcal{J}(., a, \omega)$. We have, $0 \in \partial \mathcal{J}(\Hat{f})$. Therefore, by proposition \ref{prop2.7} $\exists v\in \text{sgn}(\xi_{\Hat{f}})$ which satisfies,
\begin{align}\label{eq922}
   0 = \|\xi_{\Hat{f}}\|_{L^1(q)}\int_q v\xi_h dxdt + \epsilon\sum_{i=1}^{k}\frac{\Hat{f}_i h_i}{\|\Hat{f}\|} -\sum_{i=1}^{k}\omega_i h_i, \hspace{.5cm}  \forall h \in \mathbb{R}^k.
\end{align}
Testing \eqref{eq96a} with $\xi_h$, we get, 
\begin{align*}
    \int_{Q} U \xi_h = \sum_{i=1}^{k}\psi(x_i, T) h_i.
\end{align*}
If $\|\omega\| > \epsilon$, we choose a control $U = \|\xi_{\Hat{f}}\|_{L^1(q)}v\chi_q$, $v\in \text{sgn}(\xi_{\Hat{f}})$.
 Comparing with \eqref{eq922} we get,
 \begin{align*}
     \sum_{i=1}^{k} h_i(\psi(x_i, T) - \omega_i) &= -\epsilon\sum_{i=1}^{k}\frac{\Hat{f}_ih_i}{\|\Hat{f}\|}. 
 \end{align*}
Since the above equality holds for all $h \in \mathbb{R}^k$,
we have,
\begin{equation}
    \sum_{i=1}^{k}(\psi(x_i,T)-\omega_i)^2 \leq \epsilon. 
\end{equation}
hence the approximate controllability condition, \eqref{eq92} follows.
\end{proof}

\begin{coro}\label{cor1.7.1}
The linearised system \eqref{eq93a}-\eqref{eq93c} is approximately controllable.    
\end{coro}
\begin{proof}
We decompose the system \eqref{eq93a}-\eqref{eq93c} into the following two systems.
\begin{align}
    \psi_1' + (\Delta)^2 \psi_1 + \Delta(g(z)\psi_1) &= 0, \text{ in } \Omega \times(0,T),\label{eq913a}\\
    \frac{\partial \psi_1}{\partial \textbf{n}} = \frac{\partial \Delta\psi_1}{\partial \textbf{n}} &= 0, \text{ on } \partial\Omega \times(0,T),\label{eq913b}\\
    \psi_1(.,0) &= \psi_0, \text{ in } \Omega .\label{eq913c}
\end{align}
\begin{align}
    \psi_2' + (\Delta)^2 \psi_2 + \Delta(g(z)\psi_2) &= U\chi_q, \text{ in } \Omega \times(0,T),\label{eq914a}\\
    \frac{\partial \psi_2}{\partial \textbf{n}} = \frac{\partial \Delta\psi_2}{\partial \textbf{n}} &= 0, \text{ on } \partial\Omega \times(0,T),\label{eq914b}\\
    \psi_2(.,0) &= 0, \text{ in } \Omega.\label{eq914c}
\end{align}
The existence results of the above systems follows from the theorem \ref{thm1.1}.
  The system \eqref{eq914a}- \eqref{eq914c} is approximately controllable by theorem \ref{thm1.7}. Given $\epsilon > 0$ and state $\omega \in \mathbb{R}^k$, we modify the desired state to be $(\omega_i - \psi_1(x_i, T))_{1\leq i \leq k}$. Then there exists a control, $U \in L^2(q)$ such that,
\begin{equation*}
   \sum_{i=1}^{k} (\psi_2(x_i,T)- (\omega_i - \psi_1(x_i, T)))^2 \leq \epsilon. 
\end{equation*}  
Subsequently, $\sum_{i=1}^{k}(\psi(x_i,T)- \omega_i)^2 \leq \epsilon$.
More precisely, if the desired state, $\sum_{i=1}^{k}(\omega_i - \psi_1(x_i, T))^2 > \epsilon$, we choose a control defined as in theorem \ref{thm1.7}. Otherwise we choose the null control.  
\end{proof}


\subsection{The nonlinear problem}
In this section, we extend the results proved in the previous section to the semilinear system using a fixed point argument. More precisely, we use the Kakutani's fixed point theorem \cite{AUBIN} to prove the pointwise approximate controllability of the nonlinear problem. The proof is similar to the proof in \cite[Theorem 1]{DRA}. 


\begin{lemma}[Kakutani's fixed point theorem]
    Let $K$ be a non-empty compact convex subset of a Hilbert space $X$. Assume $\Lambda$ is an upper hemi-continuous set valued map from $K\rightarrow \mathcal{P}(K)$ with non-empty, closed convex images. Then $\Lambda$ has a fixed point. i.e., $\exists z\in K$ such that $ z\in \Lambda(z)$.
\end{lemma}

 \begin{theorem}
     Let $\varphi_0 \in H^2(\Omega)$ and the potential satisfies \emph{\textbf{[A1]}}. Then there exists a control such that the solution of the nonlinear system \eqref{eq91a}- \eqref{eq91c} satisfies the approximate controllability condition, \eqref{eq092}.
\end{theorem}

\begin{proof}   
Consider the linearised system \eqref{eq93a}-\eqref{eq93c}. For a fixed $z\in L^2(Q)$, we have by corollary \ref{cor1.7.1}, the solution of the system \eqref{eq93a}-\eqref{eq93c} with a control $U = \|\xi_{\Hat{f}}\|_{L^1(q)}v\chi_q$, $\Hat{f} = \mathcal{N}(g(z), (\omega_i - \psi_1(x_i, T))_{1\leq i \leq k})$ satisfies the approximate controllability condition,
\begin{equation}\label{eq2.350}
   \sum_{i=1}^{k} (\psi(x_i,T)- (\omega_i - \psi_1(x_i, T)))^2 \leq \epsilon. 
\end{equation} 
 Here, $\psi_1$ represents the solution of the system \eqref{eq913a}- \eqref{eq913c}. Observe that from \eqref{eq2.18},  
\begin{equation*}
    \|\xi_{\Hat{f}}\|_{L^1(q)} \leq C\|\Hat{f}\|,  
\end{equation*}
for a constant $C>0$ which depends on $\Omega$.

 The set of desired data, $$\{ (\omega_i- \psi_1(x_i, T))_{1\leq i\leq k} : \psi_1 \text{ solves the system } \eqref{eq913a} - \eqref{eq913c} \},$$ is bounded in $\mathbb{R}^k$. 
This follows from two important observations: By assumption \textbf{[A1]}, $\|g(z)\|_{L^\infty} \leq C, \forall z\in L^2(Q)$ where the constant $C$ is independent of $z$. In addition, using the energy estimate \eqref{eq2.9}, $\{\psi_1(x_i, T)\} \subseteq \mathbb{R}^k$, is bounded by a constant which depends on $\psi_0$, $\Omega$.
\\Note that, for a bounded set $B\subseteq L^\infty(Q)$ and a compact set $K \subseteq \mathbb{R}^k$ we have, $\mathcal{N}(B\times K)$ is bounded using the uniformity of the constant in the coercivity condition \eqref{eq93.22}.
Therefore, we get a uniform bound for $\|\Hat{f}\|$ and hence, 

 \begin{equation}\label{eq2.35}
    \|U\|_{L^\infty(Q)} = \|\xi_{\Hat{f}}\|_{L^1(q)} \leq C^*,
\end{equation}
for a constant $C^*>0$ independent of $z$. Let $$\mathcal{U} :=\{ U\in L^\infty(Q) : \|U\|_{L^\infty(Q)} \leq C^* \},$$
Where $C^* >0$ is the $z$ independent constant we have obtained in the previous step \eqref{eq2.35}.
Now we define a map, $\Lambda : L^2(Q) \rightarrow L^2(Q)$ as follows. For, $z\in L^2(Q)$, we define,
\begin{align}\label{eq1.33}
\Lambda(z) = \{ \psi : \psi \text{ solves } \eqref{eq93a}-\eqref{eq93c} \text{ for some } U\in \mathcal{U} \text{ and } \psi \text{ satisfies } \eqref{eq2.350}\}.
\end{align} 
 Observe that by the above construction, for every $z\in L^2(Q)$, $\Lambda(z)$ is non-empty. In order to apply Kakutani's fixed point theorem we need to verify the following hypothesis.


\begin{enumerate}[label=\roman*)]
    \item There exists a non-empty compact convex set $K\subseteq L^2(Q)$, such that $\Lambda(z) \subseteq K$ for all $z\in L^2(Q)$.
    \item For each $z\in K$, $\Lambda(z)$ is a non-empty, closed, convex subset of $K$. 
    \item The map, $\Lambda : K\rightarrow K$ is upper hemi-continuous.
\end{enumerate}

We prove all the above hypotheses one by one. 

\begin{enumerate}[label=\roman*)]
\item Let 
$$ \mathcal{S} = \{ \psi: \psi \text{ solves } \eqref{eq93a}-\eqref{eq93c} \text{ for } U\in \mathcal{U} \text{ and } z\in L^2(Q)\}.$$
Then $\mathcal{S}$ is a bounded subset of $L^2(Q)$ by applying the uniform boundedness of $g(z)$ and \eqref{eq2.35} in \eqref{eq2.9}. Now using the compact embedding of the solution space $\mathcal{V}$ defined in \ref{thm1.1} into $L^2(Q)$, we get $\mathcal{S}$ is relatively compact in $L^2(Q)$, hence there exists a compact set $K'$ such that $\forall z\in L^2(Q), \,\,\Lambda(z) \subseteq K'$. Finally we define $K$ as the convex hull of $K'$ to ensure the convexity.

\item It follows from corollary \ref{cor1.7.1} that for each $z\in K$, $\Lambda(z) \neq \emptyset$. Further, the convexity of $\Lambda(z)$ follows by  linearity of the system, convexity of $B(\omega, \epsilon)$ and $\mathcal{U}$. Since $\Lambda(z)\subseteq K$, for $K$ compact, it is enough to prove $\Lambda(z)$ is closed. Let $\psi_n \in \Lambda(z)$ be a sequence, corresponding to $U_n \in \mathcal{U}$. Then from estimates, \eqref{eq2.35}, \eqref{eq2.9} using compactness results, there exists subsequences satisfying, 
\begin{align*}
    \psi_n &\rightharpoonup \psi  \text{ in } L^2(Q),\\
    U_n &\overset{\ast}{\rightharpoonup} U  \text{ in } L^\infty(Q).
\end{align*}
By applying above convergence results, we will be able to pass to the limit for $(\psi_n, U_n)$ in the approximated problem analogous to \eqref{eq93a}-\eqref{eq93c}. This yields, $\psi$ solves the system \eqref{eq93a}-\eqref{eq93c} for control $U$. Further, consider the following system satisfied by $\rho_n:= (\psi_n-\psi)$,
\begin{align*}
    \rho_n' + (\Delta)^2\rho_n - \Delta(g(z)\rho_n) &= (U-U_n)\chi_{q}, \text{ in } \Omega \times(0,T),\\
    \frac{\partial\rho_n}{\partial\textbf{n}} = \frac{\partial\Delta\rho_n}{\partial\textbf{n}} &= 0, \text{ on } \partial\Omega \times(0,T),\\ 
    \rho_n(0) &= 0, \text{ in } \Omega.
\end{align*}
By testing with $\rho_n$, $\Delta \rho_n$, we obtain, an energy estimate analogous to \eqref{eq2.9} hold for $\rho_n$, $n \in \mathbb{N}$. Moreover, we have by interpolation embedding, $L^2(0,T;H^4(\Omega))\cap H^1(0,T;H) \hookrightarrow C([0,T];H^2(\Omega)) \hookrightarrow C([0,T];C(\Bar{\Omega}))$. Since the energy estimate hold for all $n\in \mathbb{N}$, by passing to the limit $n\rightarrow \infty$, we have, $\|U_n-U\| \rightarrow 0$, consequently, $|\rho_n(x_i,T)| \rightarrow 0$. 
Therefore we can pass to the limit in the following approximate controllability condition for $\psi_n$.
\begin{align*}
    \sum_{i=1}^{k}(\psi_n(x_i,T)- (\omega_i- \psi_1(x_i, T)))^2 &\leq \epsilon, \hspace{.5cm} \forall n\in \mathbb{N}. 
\end{align*}
 This leads to, $\sum_{i=1}^{k}(\psi(x_i,T) -(\omega_i- \psi_1(x_i, T)))^2 \leq \epsilon$. Hence, $\psi \in \Lambda(z)$ and $\Lambda(z)$ is closed.

\item A map $\Lambda$ is said to be upper hemi-continuous at $z_0$ if,
\begin{align*}
    \lim_{z_n\rightarrow z_0}\sup \sigma(\Lambda(z_n), p) \leq \sigma(\Lambda(z_0), p), \hspace{.25cm}\forall p\in L^2(Q),
\end{align*}
where $\sigma(\Lambda(z), p) := \sup_{y\in \Lambda(z)} (y, p)$. For our problem, $\Lambda(z_n)$ is compact, $\exists \psi_n \in \Lambda(z_n)$ such that, $\sigma(\Lambda(z_n), p) = (\psi_n, p)$. Observe that since $\psi_n \in K$ is  compact, there exists a limit $\psi \in K$ such that $\lim_{n\rightarrow\infty} (\psi_n, p) = (\psi, p)$. 
We can repeat the arguments in \cite{Fabre95, DRA} to conclude that $\psi \in \Lambda(z_0)$. 
 It follows immediately that,
\begin{align*}
  \lim_{z_n\rightarrow z_0}\sup \sigma(\Lambda(z_n), p) = \lim_{n\rightarrow \infty}(\psi_n, p) = (\psi, p) \leq \sup_{\Bar{\psi}\in \Lambda(z_0)}(\Bar{\psi}, p) = \sigma(\Lambda(z_0), p). 
\end{align*}
Hence the map $\Lambda$ is upper hemi-continuous.
\end{enumerate}

Now by applying Kakutani's fixed point theorem, we get a fixed point of the map $\Lambda$. It follows that,  $\exists U \in \mathcal{U}$ such that solution of \eqref{eq91a}- \eqref{eq91c} corresponding to $U$ satisfies  \eqref{eq092}. Hence the approximate controllability of the nonlinear system follows.
\end{proof}

\section{Approximate controllability with a pointwise control}

In this section, we consider a problem involving  pointwise control, along with pointwise observations. The governing PDE remains the same, however the requirement of a strong solution for a low regular control introduces additional difficulties: therefore, we restrict ourselves to the linear system. The problem is defined as follows.

\textbf{Approximate Controllability Problem Description: }Let $\{y_j : 1\leq j\leq m\}$ and $\{x_i : 1\leq i\leq k\}$ be finite disjoint subsets of $\Omega$. Given $\epsilon >0$ and a desired state $\omega = (\omega_1,\omega_2, ..., \omega_k) \in \mathbb{R}^k$, our aim is to find a control of the form, $U(x,t) = \sum_{j=1}^{m} U_j(t) \delta_{y_j}$, where $U_j \in L^r(0,T)$ for $1\leq j\leq m$, $1\leq r <\infty$ such that the solution of the system, 
\begin{align}
    \psi' + (\Delta)^2\psi - \Delta(a(x,t)\psi) &= U, \text{ in } \Omega \times(0,T),\label{eq2.1}\\
    \frac{\partial\psi}{\partial\textbf{n}} = \frac{\partial\Delta\psi}{\partial\textbf{n}} &= 0, \text{ on } \partial\Omega \times(0,T),\label{eq2.2}\\ 
    \psi(0) &= \psi_0, \text{ in } \Omega.\label{eq2.3}
\end{align}
satisfies 
\begin{equation}\label{eq2.4}
    \sum_{i=1}^{k}(\psi(x_i,T)-\omega_i)^2 \leq \epsilon. 
\end{equation}
\subsection{Wellposedness of the system}
The pointwise evaluation of $\psi$ in both space and time in condition \eqref{eq2.4} requires the continuity of the solution with respect to both the variables. To this end, we derive energy estimates for the solution of the system in higher regular spaces. We state below the existence of a solution of the linear system with a pointwise control. Due to the low regular source term in the equation, the solution is defined using the method of transposition and the proof relies on approximation techniques.

\begin{definition}\label{def2}
     Let $\psi_0 \in H$. We say $\psi \in L^{r}(0,T; H)$ solves the system \eqref{eq2.1}-\eqref{eq2.3} with $U(x,t) = \sum_{j=1}^{m} U_j(t) \delta_{y_j}$ for $U_j(t) \in L^r(0,T)$, $1\le r < \infty$, if it satisfies the following variational equation for all $f \in L^{r'}(0,T;H)$.
\begin{align}\label{eq3.5}
    \int\limits_{0}^{T}\int_\Omega \psi f dxdt - \int_\Omega \psi_0 \xi(0) dx  &= \sum_{j=1}^{m} \int\limits_{0}^{T}{}_{\mathcal{M}(\Omega)}{\langle U_j(t) \delta_{y_j}, \xi_{f} \rangle}_{C(\Bar{\Omega})} dt,
\end{align}
where $\xi_f$ represents the solution of the following adjoint system corresponding to the source term $f$. 
\begin{align}
    -\xi' + (\Delta)^2\xi - a(x,t)\Delta\xi &= f, \text{ in } \Omega \times(0,T), \label{eq2.5} \\
    \frac{\partial\xi}{\partial\textbf{n}}= \frac{\partial\Delta\xi}{\partial\textbf{n}} &= 0, \text{ on } \partial\Omega \times(0,T), \label{eq2.6}\\ 
    \xi(.,T) &= 0, \text{ in } \Omega.\label{eq2.7}
\end{align}
\end{definition}

Note that, by usual regularity argument for parabolic equation, we can show that the solution of the adjoint system with source term, $f\in L^{r'}(0,T;H)$ has the regularity, $$ \xi \in C([0,T];(H, H^2)_{1/r,r'}) \cap L^{r'}(0,T;H^2(\Omega)) \cap W^{1, r'}(0,T; H).$$ Now using the Sobolev embedding $L^{r'}(0,T;H^2(\Omega)) \hookrightarrow L^{r'}(0,T; C(\Bar{\Omega}))$ and $(H, H^2)_{1/r,r'} \hookrightarrow H$, we have \eqref{eq3.5} is well-defined. The proof follows by similar arguments as in theorem \ref{thm1.1} by Galerkin approximation and apriori estimates. Below we state and prove the existence result for the system \eqref{eq2.1}- \eqref{eq2.3} with a measure source term.
 
\begin{theorem}[Existence of solution of the  system]\label{thm3.1}
    There exists a solution of the linear system \eqref{eq2.1}-\eqref{eq2.3} corresponding to $U= \sum_{j=1}^{m}U_j(t)\delta_{y_j}$ for $U_j(t)\in L^r(0,T)$ in the sense of definition \ref{def2}.
   
   In addition, if $\psi_0 \in (W_N^{-2,p}, W_N^{2,p})_{1-\frac{1}{r}, r}$ and $a\in L^\infty(Q)$, then the solution has the following additional regularity.
    \begin{equation*}
        \psi \in L^r(0,T; W_N^{2,p}(\Omega)) \cap W^{1,r}(0,T; W_N^{-2,p}(\Omega)),
    \end{equation*}
  where  $1\leq r < \infty$, $2 < p < \infty$.
   Moreover, the solution satisfies an energy estimate,
    \begin{equation}\label{eq3.9}
        \|\psi\|_{W^{1,r}(0,T; W_N^{-2,p}(\Omega)) \cap L^r(0,T; W_N^{2,p}(\Omega))} \leq C(\|U\|_{L^r(0,T; \mathcal{M}(\Omega))} + \|\psi_{0}\|_{(W_N^{-2,p}, W_N^{2,p})_{1-\frac{1}{r}, r}}).
    \end{equation}
\end{theorem}

\begin{proof}

The first part of the theorem follows by similar arguments as in the theorem \ref{theorem1.2}. We approximate the control, $U(x,t) = \sum_{j=1}^{m} U_j(t) \delta_{y_j} \in L^r(0,T;\mathcal{M}(\Omega))$ by functions $U_{k} \in L^r(0,T; C(\Bar{\Omega}))$ such that $U_k \overset{\ast}{\rightharpoonup} U$ in $L^r(0,T;\mathcal{M}(\Omega))$ and $$\|U_k\|_{L^r(L^1(\Omega))} \leq \|U\|_{L^r(\mathcal{M}(\Omega))}.$$
Let $\psi_k$ be the solution of the corresponding linear system. Now we derive uniform priori estimates for $\psi_k$. We can get a solution of the system with pointwise control $U$ as a limit of the solution of the approximated problem. Therefore, we can conclude that there exists a transposition solution $\psi\in L^r(0,T;H)$ of the system \eqref{eq2.1}-\eqref{eq2.3} corresponding to the source term, $U(x,t) = \sum_{j=1}^{m} U_j(t) \delta_{y_j}$. 

In order to deduce the regularity, arguing as in \cite{CCK, CCK2},
consider functions $v_\beta\in \mathcal{D}(\Omega)$, $0\leq |\beta| \leq 2$ and $f \in L^{r'}(0,T; W^{-2,p'}(\Omega))$ defined as follows. For $u\in L^r(0,T;W^{2,p}(\Omega))$, 
\begin{align*}
f(u) := \sum\limits_{0\leq |\beta| \leq 2} \int\limits_{0}^{T}(D^\beta u, v_\beta).
\end{align*}
This choice of $f$ is motivated by the theorem \cite[Theorem 3.9]{ADM}.
\begin{align}
    -\xi' + (\Delta)^2\xi - a(x,t)\Delta\xi &= \sum\limits_{0\leq |\beta| \leq 2} (-1)^{|\beta|}D^\beta(v_\beta), \label{eq3.10}\\
    \frac{\partial\xi}{\partial\textbf{n}}= \frac{\partial\Delta\xi}{\partial\textbf{n}} &= 0,\\ 
    \xi(.,T) &= 0.
\end{align}
Therefore by parabolic and elliptic regularity \cite{LIO2012}, we get, $$\xi\in L^{r'}(0,T; W_N^{2,p'}(\Omega)) \cap W^{1,r'}(0,T; W_N^{-2,p'}(\Omega)) \hookrightarrow C([0,T]; (W^{-2,p'}(\Omega), W^{2,p'}(\Omega))_{1-\frac{1}{r'},r'}).$$

It is important to note that $\psi_0 \in (W_N^{-2,p}(\Omega), W_N^{2,p}(\Omega))_{1-\frac{1}{r}, r} \hookrightarrow B^{2-\frac{4}{r}}_{p,r}$ by embedding of interpolation and Besov spaces (See remark \ref{remark2}). Similarly, $\xi(0) \in (W^{-2,p'}(\Omega), W^{2,p'}(\Omega))_{1-\frac{1}{r'},r'} \hookrightarrow B^{2-\frac{4}{r'}}_{p',r'}$. Note that $B^{2-\frac{4}{r'}}_{p',r'} = B^{-(2-\frac{4}{r})}_{p',r'}$ is dual of the space $B^{2-\frac{4}{r}}_{p,r}$, therefore, the term, $\langle \psi_0, \xi_0\rangle$ is interpreted as the duality action between the spaces. 

Further, using Sobolev embedding $W^{2,p'}(\Omega) \hookrightarrow C(\Bar{\Omega})$, $1<p'$ and interpolation embedding we obtain,
\begin{equation}\label{eq3.13}
    \|\xi\|_{L^{r'}(0,T;W^{2,p'}(\Omega)) \cap C([0,T]; (W^{-2,p'}(\Omega), W^{2,p'}(\Omega))_{1-\frac{1}{r'},r'})} \leq C\sum_{0\leq |\beta| \leq 2}\|v_\beta\|_{L^{r'}(0,T;L^{p'}(\Omega))}.
\end{equation}
Now we test the adjoint equation \eqref{eq3.10} with $\psi_k$. This results, 
\begin{align*}
    \int\limits_{0}^{T} \langle f, \psi_k \rangle  dt  &= \langle \psi_0, \xi_0\rangle + \int\limits_{0}^{T}\int_\Omega U_k \xi dxdt. 
\end{align*}
Observe that by applying H\"olders to estimate the R.H.S term in the above variational equation, followed by applying the energy estimate \eqref{eq3.13}, we obtain,
\begin{align*}
    \Big|\int\limits_{0}^{T}\langle f, \psi_k \rangle \Big| = \sum_{0\leq |\beta| \leq 2}\int\limits_{0}^{T}\int_\Omega D^\beta\psi_k v_\beta dxdt  &\leq \|U_k\|_{L^r(0,T;L^1(\Omega))}\|\xi\|_{L^{r'}(0,T; C(\Bar{\Omega}))} + \|\psi_0\|_{B^{2-\frac{4}{r}}_{p,r}}\|\xi(0)\|_{B^{-(2-\frac{4}{r})}_{p',r'}}\\ &\leq C\Big(\|U\|_{L^r(0,T;\mathcal{M}(\Omega))} + \|\psi_0\|_{(W_N^{-2,p}, W_N^{2,p})_{1-\frac{1}{r}, r}} \Big)\sum_{0\leq |\beta| \leq 2}\|v_\beta\|_{L^{r'}(0,T; L^{p'}(\Omega))}. 
\end{align*}
It follows that,
\begin{equation*}
    \|\psi_k\|_{L^r(0,T;W_{N}^{2,p}(\Omega))} \leq C\Big(\|U\|_{L^r(0,T;\mathcal{M}(\Omega))} + \|\psi_0\|_{(W_N^{-2,p}, W_N^{2,p})_{1-\frac{1}{r}, r}}\Big).
\end{equation*}
Therefore, $\psi_k$ converges weakly to $\psi$ in $L^r(0,T; W^{2,p}_{N}(\Omega))$.

Using the method of transposition we can not derive the regularity in time. Therefore, we consider the following equation obtained by testing \eqref{eq2.1} with a test function $\rho \in L^{r'}(0,T; W_{N}^{2,p'}(\Omega))$,
\begin{align*}
    \langle \psi', \rho \rangle + (\Delta\psi, \Delta\rho) + (\nabla(a(x,t)\psi), \nabla\rho) = \langle U, \rho \rangle.
\end{align*}
Now estimating each term in the above equation by applying H\"olders, Gagliardo-Nirenberg and Sobolev inequalities we obtain,
\begin{align*}
    |(\Delta\psi, \Delta\rho)| &\leq \|\Delta\psi\|_{L^p} \|\Delta\rho\|_{L^{p'}} \leq \|\psi\|_{W^{2,p}_N(\Omega)}\|\rho\|_{W_N^{2,{p'}}(\Omega)},\\
    |(\nabla(a(x,t)\psi), \nabla\rho)| &= |(a(x,t)\psi, \Delta\rho)| \leq \|a\|_{L^\infty}\|\psi\|_{L^p}\|\Delta\rho\|_{L^{p'}}\\
                    &\leq \|a\|_{L^\infty}\|\psi\|_{L^p}\|\rho\|_{W_N^{2,{p'}}(\Omega)},\\
    |\langle U, \rho \rangle| &\leq \|U\|_{W_N^{-2,p}(\Omega)}\|\rho\|_{W_N^{2,{p'}}(\Omega)}.
\end{align*}
Observe that $U \in L^r(0,T;\mathcal{M}(\Omega)) \hookrightarrow L^r(0,T; W^{-2,p}(\Omega))$. Combining all the above estimates we obtain,
\begin{align*}
    \int_{0}^{T}|\langle \psi', \rho \rangle | &\leq C(\|\psi\|_{L^r(0,T; W_N^{2,p}(\Omega))}+ \|U\|_{L^r(0,T; W_N^{-2,p}(\Omega))} ) \|\rho\|_{L^{r'}(0,T; W_N^{2,{p'}}(\Omega))}.
\end{align*}
It follows that, $\psi' \in L^r(0,T; W_N^{-2,p}(\Omega))$.
\end{proof}

\begin{remark}\label{remark2}
    Let $1\leq r < \infty$, $2 < p < \infty$ be such that $2- \frac{4}{r}- \frac{2}{p} > 0$ holds. We remark that by exploiting the regularity,  $\psi \in W^{1,r}(0,T; W_{N}^{-2,p}(\Omega)) \cap L^r(0,T; W_N^{2,p}(\Omega))$ and the embedding, 
    \begin{align*}
        W^{1,r}&(0,T; W_{N}^{-2,p}(\Omega)) \cap L^r(0,T; W_N^{2,p}(\Omega)) \hookrightarrow C([0,T]; (W_{N}^{-2,p'}(\Omega), W_N^{2,p}(\Omega))_{1-\frac{1}{r}, r}) \\ &\hookrightarrow C([0,T]; B_{p,r}^{2-\frac{4}{r}}(\Omega)) \hookrightarrow C([0,T]; C_b(\Omega)),
    \end{align*}
we get a solution $\psi \in C([0,T]; C_b(\Omega))$. In the above embeddings, the first and second embedding follows from \cite[Chapter \Romannum{3}, theorem 4.10.2]{AMN},\cite{HEL}. The last embedding for the Besov space follows from \cite[Theorem 7.34]{ADM} under the conditions on $r$ and $p$ stated above.  
\end{remark}

\begin{remark}\label{remark3.3}
    Note that if $U \in L^r(0,T; W^{-2,p}(\Omega))$, then the variational formulation \eqref{eq3.5} is modified accordingly. We say $\psi \in L^r(0,T; W^{2,p}(\Omega))$ is a transposition solution if $\forall f \in L^{r'}(0,T;W^{-2,p'}(\Omega))$ the following identity is satisfied. 
\begin{align}
    \int\limits_{0}^{T} {}_{W^{2,p}(\Omega)}\langle \psi, f \rangle_{W^{-2,p'}} - \langle \psi_0, \xi(0) \rangle  &=  \int\limits_{0}^{T}{}_{W^{-2,p}}{\langle U, \xi_{f} \rangle}_{W^{2,p'}} dt.
\end{align}
As discussed earlier, $\langle \psi_0, \xi(0) \rangle$ denotes the duality action between $B^{2-\frac{4}{r}}_{p,r}$ and $B^{-(2-\frac{4}{r})}_{p',r'}$. Repeating the arguments used in theorem \ref{thm3.1}, we conclude that,
$$\psi \in L^r(0,T; W_N^{2,p}(\Omega)) \cap W^{1,r}(0,T; W_N^{-2,p}(\Omega)).$$
Moreover, the solution $\psi$ satisfies the following energy estimate.
\begin{equation}\label{eq3.15}
        \|\psi\|_{L^r(0,T; W_N^{2,p}(\Omega)) \cap W^{1,r}(0,T; W_N^{-2,p}(\Omega))} \leq C(\|U\|_{L^r(0,T; W^{-2,p}(\Omega))} + \|\psi_{0}\|_{(W_N^{-2,p}, W_N^{2,p})_{1-\frac{1}{r}, r}}).
    \end{equation}
\end{remark}

In the next subsection we prove the approximate controllability problem for  \eqref{eq2.1}-\eqref{eq2.3} with Dirac delta controls. For the proof, we will use an approximation of the measure theoretic control by $L^2$ controls and pass to the limit to devise the control. 

\subsection{The approximate controllability problem}
First we sketch the idea of the proof.
Our aim is to find a control $U$ of the form $U = \sum_{j=1}^{m} U_j(t) \delta_{y_j}$ as a limit of $L^r(0,T;H)$ functions $U_\alpha$,  such that $U_\alpha \overset{\ast}{\rightharpoonup} U$ in $L^r(0,T;\mathcal{M}(\Omega))$. For simplicity we will consider control being applied at a single point say $y_0$.
The results can be easily extended for any finite number of points. From now on we consider control, $U = U(t)\delta_{y_0}$, for $U\in L^r(0,T)$. We want to construct  $U_{\alpha} \overset{\ast}{\rightharpoonup} U$ in $L^r(0,T;\mathcal{M}(\Omega))$ such that given  $\epsilon > 0$ and $\omega = (\omega_1, \omega_2, ...,\omega_k) \in \mathbb{R}^k$,  the solution of the approximate system \eqref{eq2.1}-\eqref{eq2.3} corresponding to control $U_\alpha$ say $\psi_\alpha$ satisfies,
\begin{equation}\label{eq21}
    \sum_{i=1}^{k}(\psi_\alpha(x_i,T)-\omega_i)^2 \leq \epsilon. 
\end{equation}
Using the continuous dependence of the solution on source term \eqref{eq3.9} and by remark \ref{remark2}, we have, 
\begin{equation}\label{eq22}
  \|\psi_\alpha -\psi\|_{C([0,T]);C_b(\Omega)} \leq C\| U_\alpha -U\|_{L^r(0,T;\mathcal{M}(\Omega))}, 
\end{equation}
where $\psi$ denotes the solution of the system corresponding to control $U$. Then, by passing to the limit as $\alpha \rightarrow 0$, it follows that the pointwise control $U$ solves the pointwise approximate controllability problem. Here we assume the initial data, $\psi_0 = 0$, since for non-zero initial data the results follow immediately as in the corollary \ref{cor1.7.1}.

If we consider the controls obtained as in section 2 with supports contained in $B(0, \alpha) \times (0,T)$, the control will converge to zero in the sense of distributions as $\alpha$ tends to $0$. But, for a desired state $\omega$ satisfying  $\|\omega\| > \epsilon$, the approximate controllability can not be established using the null control.
Therefore we modify the cost functional $\mathcal{J}$ given in \eqref{eq94}. We define, $\mathcal{J}_0 : \mathbb{R}^k \rightarrow \mathbb{R}$ by, for $f =(f_1,...,f_k) \in \mathbb{R}^k$,
\begin{equation}\label{eq24}
    \mathcal{J}_0(f) = \mathcal{J}_0(f, a, \omega) := \frac{1}{2}\Big(\frac{1}{|q|}\int_q |\xi_f|dxdt\Big)^2 + \epsilon\big(\sum_{i=1}^{k}(f_i)^2\big)^{1/2} - \sum_{i=1}^{k}\omega_i f_i,
\end{equation}
where $\epsilon > 0$ and $\omega \in \mathbb{R}^k$ are given as in the statement of the control problem and $\xi_f$ represents the solution of the associated backward problem \eqref{eq95a}-\eqref{eq95c}. Observe that, the only difference in $\mathcal{J}_0$ from $\mathcal{J}$ is in a constant that appears in the first term, which is independent of $f$; we can follow the arguments in section 2 to prove that the cost $\mathcal{J}_0$ is continuous, strictly convex and coercive. Hence a minimiser exists. Further, the minimiser of $\mathcal{J}_0$, say, $\Hat{f}$ can be characterised as follows. There exists $v\in \text{sgn}(\xi_{\Hat{f}})$ such that,
\begin{align}\label{eq25}
   0 = \frac{1}{|q|^2}\|\xi_{\Hat{f}}\|_{L^1(q)}\int_q v\xi_h dxdt + \epsilon\sum_{i=1}^{k}\frac{\Hat{f}_ih_i}{\|\Hat{f}\|} -\sum_{i=1}^{k}\omega_i h_i, \hspace{.5cm}  \forall h \in \mathbb{R}^k,
\end{align} 
where $\xi_{\Hat{f}}, \xi_h$ denotes the solution of the adjoint system \eqref{eq95a}-\eqref{eq95c} with terminal data $\sum_{i=1}^{k}\Hat{f}_i\delta_{x_i}$ and $\sum_{i=1}^{k}h_i\delta_{x_i}$ respectively. By testing \eqref{eq2.1} corresponding to control $V \in L^r(0, T, H)$ on RHS, with $\xi_h$, we obtain, 
\begin{align}\label{eq26}
    \int_{Q}\xi_h V dxdt = \sum_{i=1}^{k} h_i\psi(x_i, T).
\end{align}
It follows by comparison that for control, $$V = \frac{1}{|q|^2}\|\xi_{\Hat{f}}\|_{L^1(q)} v(x,t)\chi_q \; \; \text{where} \;\; v\in \text{sgn}(\xi_{\Hat{f}}),$$  the corresponding solution of the system \eqref{eq2.1}-\eqref{eq2.3} satisfies the controllability condition, \eqref{eq2.4}. We use this technique to obtain a sequence of approximate controls $U_{\alpha}$.\\

\textbf{Approximate controllability for approximated problem:}
Given $\epsilon > 0$ and a desired state, $\omega \in \mathbb{R}^k$, consider the approximate controllability problem where the control is applied on a set, $B(y_0, \alpha) \times (0,T)$, $\alpha > 0$. Let $q = q_\alpha := B(y_0, \alpha)\times (0,T)$, correspondingly we get a control, $$U_\alpha(x,t) = \frac{1}{|q_\alpha|^2}\|\xi_\alpha\|_{L^1(q_\alpha)} v_\alpha(x,t) \chi_{q_\alpha},\;\; \text{
where} \;\; v_\alpha \in \text{sgn}(\xi_\alpha),$$ such that the solution of the approximate system $\psi_\alpha$ satisfies \eqref{eq21}. Let $\Hat{f}_\alpha$ be the unique minimiser of the cost $\mathcal{J}_0$ with $q = q_\alpha$, and denote it by $\mathcal{J}_0(., \alpha)$.  We denote the solution of the adjoint system \eqref{eq95a}-\eqref{eq95c} corresponding to $\Hat{f}_\alpha$ by $\xi_\alpha$.

Thanks to the uniform coercivity of $\mathcal{J}_0(.,\alpha),$ for fixed $\alpha_0 > 0$, $\Hat{f}_\alpha$, are uniformly bounded in $\mathbb{R}^k$ for $\alpha \leq \alpha_0$. Hence there exists a convergent subsequence of $\{\Hat{f}_\alpha\}_{0< \alpha \leq \alpha_0} \subseteq \mathbb{R}^k$ and a limit $\Hat{f}\in \mathbb{R}^k$. We prove that $U_\alpha$ converges to the function, 
\begin{equation}\label{CharU}
U = \frac{1}{T}\|\xi_{\Hat{f}}(y_0)\|_{L^1(0,T)}v(t)\delta_{y_0},
\end{equation}
in the distributional sense, where $\xi_{\Hat{f}}$ denotes the solution of adjoint problem associated with $\Hat{f}$ and $v$ represents the weak* limit of $v_\alpha \in L^\infty(Q)$ and $v\in \text{sgn}(\xi_{\Hat{f}})$. Moreover,  $\xi_{\Hat{f}} \in L^{r'}(0,T; C_b(\Omega))$ is continuous in the space variable for almost every $t$. Hence the control $U$ in \eqref{CharU} is well defined.

 Observe that $U_\alpha$, $0 < \alpha \leq \alpha_0$ are uniformly bounded. For, substituting $h= \Hat{f}_\alpha$ in \eqref{eq25} we get,
\begin{equation}\label{eq3.19}
    \Big(\frac{1}{|q_\alpha|}\|\xi_\alpha\|_{L^1(q_\alpha)}\Big)^2 \leq C(\epsilon + \|\omega\|)\|\Hat{f}_\alpha\| \leq C^*(\epsilon + \|\omega\|).
\end{equation}
The last inequality follows from the uniform boundedness of $\Hat{f}_\alpha$. 
 To prove the distributional convergence of controls $U_\alpha \rightarrow U$ as defined in \eqref{CharU},  consider  the action of $U_\alpha - U \in L^r(0,T;\mathcal{M}(\Omega))$ on $\phi \in C(\Bar{Q})$,

\begin{align}
    \langle U_\alpha - U, \phi \rangle &= \Bigl\langle \Big(\frac{1}{|q_\alpha|^2}\|\xi_\alpha\|_{L^1(q_\alpha)} v_\alpha\chi_{q_\alpha}\Big) - \frac{1}{T}\|\xi_{\Hat{f}}(y_0)\|_{L^1(0,T)}v\delta_{y_0} , \phi \Bigr\rangle \nonumber\\
    &\leq \Bigl\langle \Big(\frac{1}{|q_\alpha|}\|\xi_\alpha\|_{L^1(q_\alpha)} v_\alpha\frac{1}{|q_\alpha|}\chi_{q_\alpha}\Big) - \frac{1}{|q_\alpha|}\|\xi_\alpha\|_{L^1(q_\alpha)}v_\alpha\delta_{y_0} , \phi \Bigr \rangle \nonumber\\
    &\hspace{.5cm}+\Bigl\langle \frac{1}{|q_\alpha|}\|\xi_\alpha\|_{L^1(q_\alpha)}(v_\alpha - v)\delta_{y_0}  , \phi \Bigr \rangle + \Bigl\langle \frac{1}{|q_\alpha|}\|\xi_\alpha\|_{L^1(q_\alpha)}v \delta_{y_0} - \frac{1}{T}\|\xi_{\Hat{f}}(y_0)\|_{L^1(0,T)}v\delta_{y_0} , \phi \Bigr \rangle \nonumber\\
     &\leq  \Bigl\langle \frac{1}{|q_\alpha|}\|\xi_\alpha\|_{L^1(q_\alpha)} v_\alpha\Big(\frac{1}{|q_\alpha|}\chi_{q_\alpha}- \delta_{y_0}\Big) , \phi \Bigr \rangle \nonumber\\
     &\hspace{.5cm}+\Bigl\langle \frac{1}{|q_\alpha|}\|\xi_\alpha\|_{L^1(q_\alpha)}(v_\alpha - v)\delta_{y_0}  , \phi \Bigr\rangle + \Bigl\langle  \Big(\frac{1}{|q_\alpha|}\|\xi_\alpha\|_{L^1(q_\alpha)} - \frac{1}{T}\|\xi_{\Hat{f}}(y_0)\|_{L^1(0,T)} \Big) v\delta_{y_0} , \phi \Bigr \rangle.\label{eq3.20}
\end{align}
We denote by $\delta_{y_0}$ the characteristic function of the set $\{y_0\} \times [0,T]$. Let $I_1, I_2, I_3$ denote the terms on the R.H.S of the last inequality. We estimate $I_j$ for each $j$ and evaluate the limit as $\alpha$ tends to zero. Observe that,
\begin{align*}
    |I_1| &= \Big|\Bigl\langle \frac{1}{|q_\alpha|}\|\xi_\alpha\|_{L^1(q_\alpha)} v_\alpha\Big(\frac{1}{|q_\alpha|}\chi_{q_\alpha}- \delta_{y_0}\Big) , \phi \Bigr \rangle \Big|\\
    &\leq C \Big|\Bigl\langle \frac{1}{|q_\alpha|}\chi_{q_\alpha}- \delta_{y_0} , \phi \Bigr \rangle\Big|  \longrightarrow 0, 
\end{align*}
as $ \alpha\rightarrow 0$. The convergence follows by  the estimate \eqref{eq3.19}
and the distributional convergence, $\frac{1}{|q_\alpha|}\chi_{q_\alpha} \rightarrow \delta_{y_0}$. 
\begin{align*}
    |I_2| &= |\Bigl\langle \frac{1}{|q_\alpha|}\|\xi_\alpha\|_{L^1(q_\alpha)}(v_\alpha - v)\delta_{y_0}  , \phi \Bigr \rangle| \leq C|\langle (v_\alpha - v)\delta_{y_0}  , \phi  \rangle|.  
\end{align*}
Using the convergence $v_\alpha \overset{\ast}{\rightharpoonup} v$ in $L^\infty(Q)$, we have $I_2$ converges to $0$ as $\alpha\rightarrow 0$.
\begin{align*}
   |I_3| &= |\Bigl\langle  \Big(\frac{1}{|q_\alpha|}\|\xi_\alpha\|_{L^1(q_\alpha)} - \frac{1}{T}\|\xi_{\Hat{f}}(y_0)\|_{L^1(0,T)} \Big) v\delta_{y_0} , \phi \Bigr \rangle |
\end{align*}
We have $|v(x,t)|\leq 1$ a.e in $Q$. Observe that,
\begin{align}
    \Big|\frac{1}{|q_\alpha|}\|\xi_\alpha\|_{L^1(q_\alpha)} - \frac{1}{T}\|\xi_{\Hat{f}}(y_0)\|_{L^1(0,T)} \Big| &=  \frac{1}{T}\Big|\int\limits_{0}^{T}\Big(\frac{1}{|B(y_0, \alpha)|}\int_{B(y_0, \alpha)}(|\xi_\alpha(x,s)| - |\xi_{\Hat{f}}(y_0,s)|)dx\Big)ds\Big| \nonumber\\
    &\leq \frac{1}{T}\Big|\int\limits_{0}^{T}\Big(\frac{1}{|B(y_0, \alpha)|}\int_{B(y_0, \alpha)}(|\xi_\alpha(x,s)| - |\xi_{\Hat{f}}(x,s)|)dx\Big)ds\Big|\nonumber\\
    &\hspace{.5cm}+\frac{1}{T}\Big|\int\limits_{0}^{T}\Big(\frac{1}{|B(y_0, \alpha)|}\int_{B(y_0, \alpha)}|\xi_{\Hat{f}}(x,s)| - |\xi_{\Hat{f}}(y_0,s)|dx\Big)ds\Big|\nonumber\\
    &\leq \underbrace{\frac{1}{T}\Big|\int\limits_{0}^{T}\Big(\frac{1}{|B(y_0, \alpha)|}\int_{B(y_0, \alpha)}|\xi_\alpha(x,s) - \xi_{\Hat{f}}(x,s)|dx\Big)ds\Big|}_{k_1} \nonumber\\
    &\hspace{.5cm}+\underbrace{\frac{1}{T}\Big|\int\limits_{0}^{T}\Big(\frac{1}{|B(y_0, \alpha)|}\int_{B(y_0, \alpha)}|\xi_{\Hat{f}}(x,s) - \xi_{\Hat{f}}(y_0,s)|dx\Big)ds\Big|}_{k_2}\label{eq3.24}
\end{align}
Now using the continuous dependence of solution of the adjoint system on terminal data from \eqref{eq2.15}, and the embedding $W^{2,p'}(\Omega) \hookrightarrow C(\Bar{\Omega})$, we get,
\begin{equation}
\|\xi_\alpha - \xi_{\Hat{f}}\|_{L^{r'}(0,T; C(\Bar{\Omega}))} \leq C\|\Hat{f}_\alpha - \Hat{f}\|.
\end{equation}
We conclude that the term $k_1$ of \eqref{eq3.24} vanishes as $\alpha \rightarrow 0$. Observe that, $\xi_{\Hat{f}}$ is continuous in the space variable for almost every $t\in (0,T)$. Therefore using the distributional convergence $\frac{1}{|B(y_0, \alpha)|}\int_{B(y_0, \alpha)}|\xi_{\Hat{f}}(y)| dy \rightarrow |\xi_{\Hat{f}}(y_0)|$, the term $k_2$ in \eqref{eq3.24} converges to $0$. Consequently, $|I_3| \rightarrow 0$ as $\alpha \rightarrow 0$. Hence from \eqref{eq3.20}, $U_\alpha \overset{\ast}{\rightharpoonup} U$ in $L^r(0,T;\mathcal{M}(\Omega))$.

 It is left to prove that $v\in \text{sgn}(\xi_{\Hat{f}})$. We have by unique continuation theorem \ref{theorem2.4}, that if $\Hat{f} \neq 0$, then $\xi_{\Hat{f}}\neq 0$ except on a set of measure $0$. Observe that if $\xi_{\Hat{f}}(x,t) \neq 0$, there exist a sequence such that $\text{sgn}(\xi_\alpha)(x,t) = \text{sgn} (\xi_{\Hat{f}}(x,t))$ for sufficiently small $\alpha$. We get $\text{sgn}(\xi_\alpha) \rightarrow \text{sgn}(\xi_{\Hat{f}})$ for a.e. $(x,t) \in \Omega \times (0,T)$. Therefore by uniqueness of the limit it follows that $v\in \text{sgn}(\xi_{\Hat{f}})$. 
 We now summarize the above discussion in the following proposition. 
 
\begin{prop}\label{prop3.3}
    Let $\epsilon>0$, $\omega\in \mathbb{R}^k$ be given. There exists a sequence of controls $U_\alpha \in L^r(0,T;H),$ $\alpha >0$ such that the solution of the system  \eqref{eq2.1}-\eqref{eq2.3} corresponding to $U_\alpha$ satisfies the approximate controllability controllability condition \eqref{eq21}. And $U_\alpha$ converges to $U \neq 0$ in $L^r(0,T;\mathcal{M}(\Omega))$  as $\alpha \rightarrow 0$ with support of $U$ contained in $\{y_0\} \times (0,T)$. 
\end{prop}

\subsection{The approximate controllability problem with pointwise control}
We now present the main result of this section, establishing the approximate controllability of the Cahn–Hilliard system with pointwise observations and pointwise controls. The proof follows from the proposition \ref{prop3.3}.

\begin{theorem}\label{thm3.3}
 Let $1\leq r < \infty$, $2 < p < \infty$ be such that $2- \frac{4}{r}- \frac{2}{p} > 0$ holds. Given $\epsilon>0$, $\omega\in \mathbb{R}^k$. There exists a pointwise control $U = \sum_{j=1}^{m}U_j(t)\delta_{y_j} \in L^r(0,T; \mathcal{M}(\Omega))$, $U_j(t)\in L^r(0,T)$ such that the corresponding solution of \eqref{eq2.1}-\eqref{eq2.3}, $\psi \in L^r(0,T;W_N^{2,p}(\Omega)) \cap W^{1,r}(0,T; W_N^{-2,p}(\Omega))$ satisfies the approximate controllability condition \eqref{eq2.4}.
\end{theorem}

\begin{proof}
W.l.g we assume that $U$ is a control acting at a single point $y_0$. If the desired data, $\omega$ satisfies $\|\omega\| \leq \epsilon$, then \eqref{eq2.4} holds trivially with the null control $U= 0$.
For the case, $\|\omega\|>\epsilon$, we establish the existence of a pointwise control using proposition \ref{prop3.3}. The continuous dependence estimate \eqref{eq3.15} enable us to pass to the limit in \eqref{eq21} to obtain the approximate controllability condition \eqref{eq2.4}. This completes the proof.
\end{proof}

\addcontentsline{toc}{chapter}{Bibliography}
\nocite{*}
\printbibliography

\end{document}